\documentclass[preprint,number]{elsarticle}

\usepackage{geometry}
\usepackage{amsmath,amssymb,amsthm,mathtools}
\usepackage{booktabs}
\usepackage{array}
\usepackage{float}
\usepackage{caption}
\usepackage{algorithm}
\usepackage{algpseudocode}
\usepackage{xcolor}
\usepackage{hyperref}
\usepackage{bookmark}
\hypersetup{hidelinks,hypertexnames=false}
\biboptions{sort&compress}
\newcommand{\doi}[1]{\url{https://doi.org/#1}}
\newtheorem{theorem}{Theorem}[section]

\newtheorem{proposition}[theorem]{Proposition}

\newtheorem{remark}[theorem]{Remark}

\newtheorem{definition}[theorem]{Definition}

\begin{document}

\begin{frontmatter}

\title{Higher-order stabilized Runge--Kutta methods \\ suitable for solution limiting}

\author[oden,ase]{Sreeram Shankar\corref{cor1}}
\ead{sreeramshankar08@utexas.edu}
\author[oden,ase]{Jesse Chan}
\ead{jesse.chan@oden.utexas.edu}
\cortext[cor1]{Corresponding author.}
\affiliation[oden]{organization={Oden Institute for Computational Engineering and Sciences,
The University of Texas at Austin},
            city={Austin},
            state={Texas},
            country={USA}}
\affiliation[ase]{organization={Department of Aerospace Engineering and Engineering Mechanics,
The University of Texas at Austin},
            city={Austin},
            state={Texas},
            country={USA}}

\begin{abstract}
Existing higher-order stabilized Runge-Kutta time-stepping methods, such as ROCK and DUMKA, have high internal amplification and do not possess monotone-in-time stages, which makes them unusable for solution limiting. 
We introduce Chain Runge--Kutta Gegenbauer (CRKG) methods as a family of stabilized Runge-Kutta schemes suitable for solution limiting. CRKG schemes combine several orthogonal polynomial recurrences in a step to satisfy higher-order conditions. 
The resulting methods achieve orders 3 and 4 on nonlinear problems, while retaining physical, monotone-in-time stages and low internal amplification for solution limiting.
The methods are demonstrated on viscous 
compressible flow problems, where explicit methods are restricted by the parabolic step-size limit and existing higher-order stabilized methods are destabilized by stagewise solution limiting. 
\end{abstract}


\end{frontmatter}

\section{Introduction}
\label{sec:intro}

Stabilized explicit Runge--Kutta methods (also referred to as super-time-stepping
methods~\cite{alexiades1996}) have emerged as popular and efficient methods for the
solution and simulation of stiff systems of differential
equations~\cite{verwer1996,hairer1996,abdulle2015encyclopedia}.
These methods are built on an orthogonal polynomial approximation to the stability
function $R(z)$ that, for a stage count $s$, typically extends the stability region
along the negative real axis by a factor proportional to $s^{2}$~\cite{vdh1980,meyer2014}.
This extended stability region makes them suitable for systems such as parabolic
partial differential equations, whose eigenvalues lie on or near the negative real
axis and can therefore be integrated both stably and efficiently compared with traditional explicit methods.
These methods have been well studied at low order, but
higher-order constructions face limitations made precise below.

The Runge--Kutta--Chebyshev (RKC) methods~\cite{vdh1980,verwer1990,sommeijer1998} are
the classic stabilized methods and are limited to order~2.
The Orthogonal Runge--Kutta--Chebyshev (ROCK)~\cite{abdulle2001,abdulle2002} and
DUMKA~\cite{lebedev1994,medovikov1998}
methods reach order~4.
ROCK methods achieve this by appending Runge--Kutta-like finishing stages to the end
of an orthogonal polynomial recurrence, while DUMKA methods compose steps built from
pairs of roots of the stability polynomial.
Both are efficient, but they suffer from high internal amplification, which is the propagation of small errors in the internal stages to the solution~\cite{ketcheson2014,vdh1980}, and they lack physical, monotone-in-time stages in the sense of Definition~\ref{def:physical-stages} below.
The Runge--Kutta--Legendre (RKL)~\cite{meyer2014,meyer2012} and
Runge--Kutta--Gegenbauer (RKG)~\cite{skaras2021} families are built for orders~1
and~2 and possess the convex monotone property (CMP), with RKG designed to handle
problems with Dirichlet boundary conditions; however, they are restricted to order~2.
The factorized RKG (FRKG) methods of O'Sullivan~\cite{osullivan2019}, which build on
the factorized RKC methods~\cite{osullivan2015,osullivan2017}, provide arbitrary linear
order but only second-order nonlinear order.
Beyond second order, these methods require splitting into a linear and a nonlinear component.
Such a splitting involves negative or complex substeps, can lead to order reduction, and is less
applicable generally~\cite{sheng1989,suzuki1991}.
These methods also differ in that, unlike a recurrence, the roots of the orthogonal polynomial are used as the step sizes of stagewise Euler steps. However, since the roots are generally complex, the intermediate stages are complex-valued as well.

The stabilized explicit Runge--Kutta (SERK) methods of Mart\'in-Vaquero and
coauthors~\cite{martin2009} form another low-order family
that is similar to the one detailed above in properties, limited to second order.
The extrapolated SERK (ESERK) methods~\cite{martin2016,martin2019,martin2020}
reach up to order~6 by Richardson extrapolation of lower-order SERK solutions.
Since extrapolation can raise the order of any one-step method, including the
forward Euler method~\cite[Sec.~II.9]{hairer1993}, this is regarded as a separate mechanism rather
than ``native'' high order.
Furthermore, keeping the extrapolated stability function bounded
requires heavily damped SERK stability functions, which significantly reduces the
stability interval along the real axis~\cite{martin2016,martin2020frontiers}, and
higher order requires additional storage and function evaluations.
Extrapolation-based methods are therefore not considered further.

Higher temporal order is worth pursuing in this setting for accuracy and
efficiency.
For unsteady compressible flow, temporal integrators of order three to five have been found more efficient than second-order schemes at engineering accuracy~\cite{bijl2002,kennedy2000}, and with high-order spatial discretizations such as discontinuous Galerkin methods, a second-order integrator at the step sizes permitted by stability can contribute a temporal error that exceeds the spatial error~\cite{wang2013}.
Thus, the goal of this paper is to construct a family of stabilized Runge--Kutta methods that:
\begin{enumerate}
\item achieve up to fourth order ``natively'', without extrapolation or splitting;
\item minimize internal amplification, especially compared with ROCK and DUMKA, the
existing methods of this order;
\item have physical, monotone-in-time stages in the sense of Definition~\ref{def:physical-stages},
for reasons detailed below; and
\item are practical and efficient for large-scale simulation, with the viscous
Navier--Stokes equations as a particular focus.
\end{enumerate}
To our knowledge, no existing stabilized method satisfies all four criteria.
RKC, RKL, and RKG are limited to second order, ROCK and DUMKA reach fourth order at the cost of high internal amplification and without physical, monotone-in-time stages, and FRKG and ESERK rely on splitting or extrapolation.
CRKG is thus a third family of stabilized Runge--Kutta methods to reach order four natively, after ROCK and DUMKA, and is distinguished by combining that order with physical, monotone-in-time stages and low internal amplification.

To achieve the above goal, we introduce the Chain Runge--Kutta--Gegenbauer (CRKG)
methods.
Section~\ref{sec:prelim} reviews stabilized Runge--Kutta methods and discusses their
applicability to internal stage limiting.
Section~\ref{sec:crkg} extends Runge--Kutta--Gegenbauer methods beyond second order
by chaining Gegenbauer recurrences, matching a prescribed stability polynomial, and
satisfying nonlinear order conditions beyond second order, with stage sharing and refinement
for efficiency and low amplification.
Section~\ref{sec:computation} describes the step evaluation.
Section~\ref{sec:discussion} compares the stage and amplification properties against
DUMKA3 and ROCK4, and Section~\ref{sec:numerics} demonstrates the methods on
stagewise-limited problems, including the viscous shock tube of Daru and
Tenaud~\cite{daru2001}.
\ref{app:mixed} discusses a mixed-precision evaluation made practical by the low amplification, including a GPU timing study, and \ref{app:crkg5} details a potential fifth-order construction and explains why it is not pursued.

\section{Preliminaries}
\label{sec:prelim}

\subsection{Review of stabilized Runge--Kutta methods}
\label{sec:prelim-review}

We begin by reviewing relevant properties of existing second order stabilized Runge-Kutta methods, as well as current higher-order extensions.

\subsubsection{Stabilized explicit Runge--Kutta methods}
\label{sec:prelim-stabilized}

Consider an ordinary differential equation $y'=f(y)$ and a single step of size~$h$ from $y_n$ to $y_{n+1}$.
An $s$-stage explicit Runge--Kutta method produces intermediate states
$Y_0=y_n,Y_1,\ldots,Y_s$ and a step update built from $f(Y_j)$.
On a linear test equation $y'=\lambda y$ with $z=h\lambda$, each stage and the
step are polynomials in~$z$:
\begin{equation}
\label{eq:stage-poly}
  Y_j = R_j(z)\,y_n,\qquad
  y_{n+1}=R(z)\,y_n,
\end{equation}
where $R_s\equiv R$ is the stability function of the method.
The method is stable for a given $z$ when $|R(z)|\le 1$.
For parabolic problems the relevant spectrum lies far along the negative real axis, so
a design for $R$ that maximizes stability for these problems would satisfy a condition of the form
\begin{equation}
\label{eq:real-stability}
  |R(x)|\le 1
  \qquad\text{for all }x\in[-\beta,0],
\end{equation}
with the interval length~$\beta=\beta(s)$ as large as possible.
Consistency of order~$p$ at the origin requires $R(z)=e^{z}+O(z^{p+1})$ as $z\to 0$.
A method is thus called a stabilized Runge--Kutta method when $\beta(s)$ grows like $O(s^{2})$, which
appears to be the best possible scaling for an $s$-stage explicit scheme~\cite{vdh1980}.

\subsubsection{Orthogonal-polynomial constructions}
\label{sec:prelim-orthogonal}

The standard way to obtain such an~$R$ is through a sequence of orthogonal polynomials
that is already bounded on a fixed interval.
Let $\{P_j\}$ be a family of orthogonal polynomials on $[-1,1]$, typically Jacobi
polynomials such as the Chebyshev, Legendre, or Gegenbauer polynomials---recovering,
respectively, the RKC~\cite{vdh1980,verwer1990,sommeijer1998},
RKL~\cite{meyer2014,meyer2012}, and RKG~\cite{skaras2021} methods---with a
three-term recurrence
\begin{equation}
\label{eq:three-term}
  P_0\equiv 1,\qquad
  P_1(w)=a_1 w+b_1,\qquad
  P_j(w)=(a_j w+b_j)P_{j-1}(w)-g_j P_{j-2}(w).
\end{equation}
To mitigate the oscillations of these polynomials inside the interval, in particular
for the Runge--Kutta--Chebyshev methods, a so-called damping shift
$w_0=1+\varepsilon/s^{2}$ is introduced~\cite{vdh1980,sommeijer1998}.
With the affine map
\begin{equation}
\label{eq:affine-map}
  w(z)=w_0+w_1 z,\qquad
  w_1=\frac{1+w_0}{\beta},
\end{equation}
one sets
\begin{equation}
\label{eq:Rj-orthogonal}
  R_j(z)=\frac{P_j\bigl(w(z)\bigr)}{P_j(w_0)}.
\end{equation}
Then $|R_j|\le 1$ on $[-\beta,0]$ whenever $w([-\beta,0])$ lies in the interval
where $|P_j|$ is bounded, and the recurrence~\eqref{eq:three-term} translates
into a three-term Runge--Kutta recurrence for the stages~$Y_j$.
The corresponding step update is of the form~\eqref{eq:crkg-update} (with a single recurrence, $n_c=1$) for the RKG methods.
Damping $\varepsilon>0$ slightly shortens~$\beta$ but smooths $|R|$ inside the
interval~\cite{vdh1980,sommeijer1998}.

In this setting, each stage polynomial $R_j$ is real, bounded on $[-\beta,0]$, and satisfies $R_j(z)=1+c_jz+O(z^2)$ with $c_j=R_j'(0)=w_1P_j'(w_0)/P_j(w_0)\ge 0$, so each stage approximates the solution at the time $t_n+c_jh$.
Moreover, since $w_0$ lies at or beyond the largest zero of every $P_j$, the interlacing of the zeros of orthogonal polynomials~\cite{szego1975} gives $c_j>c_{j-1}$, so the stages advance monotonically in time.
These are the properties formalized in Definition~\ref{def:physical-stages}.

\subsubsection{Higher-order stabilized methods}
\label{sec:prelim-other}

The methods above are limited to order two; we now detail the higher-order
stabilized methods.

\paragraph{ROCK}
ROCK methods~\cite{abdulle2001,abdulle2002} also use an orthogonal-polynomial
recurrence, with polynomials orthogonal with respect to a weight function chosen so
that the stability function factors as the recurrence polynomial times a fixed
low-degree polynomial.
After the recurrence stages, a short Runge--Kutta finishing procedure realizes this
factor and raises the nonlinear order to four.
The recurrence stages are of the form~\eqref{eq:Rj-orthogonal}, while the finishing
stages are not.

\paragraph{DUMKA}
DUMKA methods~\cite{lebedev1994,medovikov1998} factor a high-degree stability
polynomial into linear and quadratic factors over its roots and realize each factor
by a short explicit Runge--Kutta substep, with the order of the factors chosen to
limit the growth of intermediate values.
The intermediate states are therefore partial products of these factors, rather than
the recurrence form detailed in the previous section.

\subsubsection{Factorized stabilized methods}
\label{sec:prelim-factorized}

The factorized Runge--Kutta--Gegenbauer (FRKG) methods of
O'Sullivan~\cite{osullivan2019}, building on the factorized Runge--Kutta--Chebyshev
methods~\cite{osullivan2015,osullivan2017}, form a different class.
Each stage is a forward Euler step whose step size is derived from a root of the stability polynomial, and the
steps are ordered to limit the growth of intermediate values.
Since the roots are generally complex, so are the step sizes and the intermediate
states.
The methods have arbitrary linear order but native nonlinear order two, with higher
nonlinear order requiring splitting~\cite{osullivan2017}.
Thus, by the same standard we apply to ESERK, we do not consider splitting to be a
``native'' mechanism for high-order.
As FRKG additionally does not follow the recurrence-based construction of
Section~\ref{sec:prelim-orthogonal}, we do not focus on these methods further
in this paper.

They are, however, central to our work in another respect: the RKG stability
polynomials introduced in~\cite{osullivan2019}, formed from Gegenbauer polynomials
whose degrees are multiples of a segment count, are the polynomials on which we base
the construction of CRKG in Section~\ref{sec:crkg}.
The chaining technique we introduce in the following sections
can also be applied to factorized methods, and is the subject of future work.

\subsection{Properties relevant to solution limiting}
\label{sec:prelim-applicability}

The motivation for CRKG schemes is the design of stabilized RK methods suitable for solution limiting. Higher-order discretizations of conservation laws can produce nonphysical states near shocks and under-resolved gradients, which are typically addressed using solution limiters to enforce admissibility while retaining conservation~\cite{zhang2010maximum, zhang2017positivity, liu2024optimization, liu2024simple, kucharik2003repair, shashkov2004repair, guermond2026crl, guermond2018second, guermond2021second, pazner2021sparse, rueda2022subcell}.

The following sections review properties relevant to internal stage limiting. 
Bound-preserving limiters are often designed for a single forward Euler step, with extensions to higher order time-stepping relying on strong stability preserving (SSP) Runge--Kutta
methods~\cite{shu1988efficient, gottlieb2001strong} for which each stage can be written as a convex combination of forward Euler steps. The SSP property is therefore sufficient for stagewise limiting, but is not necessary.
Mass-redistribution-based corrections such as Liu--Zhang
limiting~\cite{liu2024optimization, liu2024simple}, the repair methods of~\cite{kucharik2003repair, shashkov2004repair},
and related conservative redistribution limiting~\cite{guermond2026crl}
remove the requirement that the limiter act on a single forward Euler step, so one is not restricted to SSPRK methods or to reformulations of Runge--Kutta methods in that
spirit~\cite{ern2022invariant}. 

However, what remains necessary is that each stage is an admissible approximation reached by advancing forward in time, and that solution perturbations introduced by limiters are not amplified in the step output. The first condition is motivated by the following: a stage obtained by a forward Euler step of negative step size integrates the 
problem backward in time and reverses its dissipation. For a parabolic operator, this step is anti-diffusive. 
We note that, while negative-coefficient SSP methods can use backward-in-time downwind stepping to overcome order barriers~\cite{shu1988efficient,gottlieb2001strong}, this
approach relies on reversing spatial upwinding and does not apply to a physical parabolic term.



\subsubsection{Physical, monotone-in-time stages}
\label{sec:prelim-physical}

\begin{definition}[Physical, monotone-in-time stages]
  \label{def:physical-stages}
  Consider a stabilized explicit Runge--Kutta method in which each stage $Y_j$ is formed as a linear combination of previously computed states $Y_i$, $i\in\mathcal{S}_j$, and of right-hand side evaluations at these states, with stage polynomials $Y_j=R_j(z)\,y_n$ as in~\eqref{eq:stage-poly}.
  The method has \emph{physical, monotone-in-time stages} on $[-\beta,0]$ if, for every stage~$j$,
  \begin{enumerate}
  \item[(i)] $R_j$ has real coefficients;
  \item[(ii)] $R_j(z)=1+c_jz+O(z^2)$ with $c_j\ge 0$, so that $Y_j$ is a consistent approximation of the solution at the time $t_n+c_jh$;
  \item[(iii)] $|R_j(x)|\le 1$ for all $x\in[-\beta,0]$; and
  \item[(iv)] $c_j>c_i$ for every $i\in\mathcal{S}_j$, so that every stage approximates the solution at a later time than the states from which it is computed.
  \end{enumerate}
\end{definition}

Conditions~(i)--(iii) alone are what we refer to as physical stages, which means that each stage is a real-valued, consistent, and bounded approximation of the solution.
Condition~(iv) adds the monotone-in-time requirement, that every stage is reached only by advancing forward in time.
A stage computed by stepping backward from a later state can still be physical in this sense, so~(iv) is not implied by~(i)--(iii).
A backward step reverses the dissipation of the problem, which is incompatible with certain applications and stability requirements, and we therefore include~(iv) in the full definition.
Condition~(iv) also depends on the form in which the method is implemented, as does the internal amplification~\cite{ketcheson2014}.

Recurrence-type methods built as in Section~\ref{sec:prelim-orthogonal} satisfy Definition~\ref{def:physical-stages}, including with damping, while FRKG violates condition~(i), since its step sizes and intermediate states are generally complex.
The higher-order methods of Section~\ref{sec:prelim-other} also fail the full definition.
ROCK4 has physical stages (conditions~(i)--(iii)), as we verify in Section~\ref{sec:discussion-physical}, but violates condition~(iv). Its finishing procedure begins with a forward Euler step of negative step size, $a_{21}<0$ for every tabulated degree~\cite{abdulle2002}, so that the first finishing stage approximates the solution at an earlier time than the last recurrence stage from which it is computed.
DUMKA violates condition~(iii), with intermediate states exceeding one by several orders of magnitude (Section~\ref{sec:discussion-physical}), as well as condition~(iv), since the last substep of every degree likewise begins with a forward Euler step of negative step size~\cite{medovikov1998}.

The low-order RKL and RKG schemes~\cite{meyer2014,skaras2021} satisfy conditions~(i)--(iii) and, in fact, the stronger convex monotone property (CMP).
In particular, every stage remains in the convex hull of previous values, which under suitable conditions on the ODE implies a discrete maximum principle.
CMP is thus sufficient for the boundedness of condition~(iii), but not necessary.

\subsubsection{Internal amplification}
\label{sec:prelim-amp}

A related but distinct notion is how strongly a perturbation of an internal
stage reaches the step output.
Following van~der~Houwen and Sommeijer~\cite{vdh1980} and
Ketcheson et al.~\cite{ketcheson2014}, consider the linear test $y'=\lambda y$ with
$z=h\lambda$ and a perturbation of size~$\epsilon$ applied to stage~$Y_j$.
Writing $y_{n+1}$ and $\widetilde y_{n+1}$ for the clean and perturbed step
outputs, the stage amplification factor is
\begin{equation}
\label{eq:theta-j}
  \theta_j(z)
  :=
  \lim_{\epsilon\to 0}
  \frac{\widetilde y_{n+1}-y_{n+1}}{\epsilon}
  =
  \frac{\partial y_{n+1}}{\partial Y_j}(z).
\end{equation}
\begin{definition}[Internal amplification]
\label{def:internal-amp}
The internal amplification of an $s$-stage method on a real stability interval
$[-\beta,0]$ is
\begin{equation}
\label{eq:internal-amp}
  \mathcal{A}
  :=
  \max_{1\le j\le s}
  \,
  \max_{\substack{z\in[-\beta,0]\\ |R(z)|\le 1}}
  \bigl|\theta_j(z)\bigr|.
\end{equation}
\end{definition}
Small $\mathcal{A}$ means that roundoff, limiter corrections, and other solution
perturbations at the stages remain small in $y_{n+1}$.
Higher-order stabilized methods fare poorly on this measure,
as will be demonstrated in Section~\ref{sec:discussion-amp}.


\section{Higher-order Chain Runge--Kutta--Gegenbauer methods}
\label{sec:crkg}

For the Runge--Kutta--Gegenbauer (RKG) methods of Skaras et
al.~\cite{skaras2021}, a single Gegenbauer three-term
recurrence yields stages satisfying Definition~\ref{def:physical-stages} and a
stability function built from Gegenbauer polynomials.
In particular, matching a prescribed stability polynomial of degree~$L$ fixes the $L$ weights
of a single recurrence uniquely, so no freedom remains for higher-order conditions beyond order two.
To obtain native nonlinear orders three and four while retaining physical,
monotone-in-time stages and a Gegenbauer-type stability polynomial, we seek
methods that combine several such recurrences in a single step.
Writing $n_c$ for the number of recurrences and $Y_{c,0}=y_n,Y_{c,1},\ldots,Y_{c,L-1}$
for the stages of the $c$th recurrence, the step update is
\begin{equation}
\label{eq:crkg-update}
  y_{n+1}=y_n+h\sum_{c=1}^{n_c}\sum_{j=0}^{L-1}b_{c,j}\,f(Y_{c,j}),
\end{equation}
with free weights $b_{c,j}$.
On the linear test equation, the stability function is
\begin{equation}
\label{eq:crkg-stability}
  R(z)=1+z\sum_{c=1}^{n_c}\sum_{j=0}^{L-1}b_{c,j}\,R_j^{(c)}(z),
\end{equation}
where $R_j^{(c)}$ is the stage polynomial of stage~$j$ in recurrence~$c$.
When $n_c=1$, both~\eqref{eq:crkg-update} and~\eqref{eq:crkg-stability} reduce
to a single Gegenbauer recurrence of the same type as the existing RKG methods.
With $n_c>1$, the additional weights supply the freedom needed for higher-order
conditions, while $R(z)$ remains a combination of Gegenbauer stage polynomials
and can be matched to a Gegenbauer stability polynomial of the desired order.
Each Gegenbauer recurrence in this formulation is called a \emph{chain}, and the resulting family is the
Chain Runge--Kutta--Gegenbauer (CRKG) methods.

\subsection{A single chain}
\label{sec:crkg-chain}

For $\nu>0$, the Gegenbauer polynomials $C_n^{(\nu)}$ are orthogonal on $[-1,1]$ with
respect to the weight $(1-w^2)^{\nu-1/2}$ and satisfy the three-term
recurrence~\cite{szego1975,dlmf}
\begin{equation}
\label{eq:gegenbauer}
  C_0^{(\nu)}(w)=1,\qquad
  C_1^{(\nu)}(w)=2\nu w,\qquad
  (n+1)\,C_{n+1}^{(\nu)}(w)=2(n+\nu)\,w\,C_n^{(\nu)}(w)-(n+2\nu-1)\,C_{n-1}^{(\nu)}(w).
\end{equation}
The parameter $\nu$ interpolates the common families: $\nu=1/2$ gives the Legendre
polynomials, $\nu=1$ the Chebyshev polynomials of the second kind, and as
$\nu\to 0$, $C_n^{(\nu)}/\nu\to(2/n)\,T_n$.
Their derivatives are again Gegenbauer~\cite{szego1975},~\cite[§18.9]{dlmf},
\begin{equation}
\label{eq:gegenbauer-derivative}
  \frac{d^{j}}{dw^{j}}\,C_n^{(\nu)}(w)=2^{j}(\nu)_j\,C_{n-j}^{(\nu+j)}(w),
\end{equation}
where $(\nu)_j=\nu(\nu+1)\cdots(\nu+j-1)$ is the Pochhammer symbol.

A \emph{chain} is a single Gegenbauer three-term recurrence of the same type as the
classical constructions in Section~\ref{sec:prelim-orthogonal}, and in particular of
the same type as RKG~\cite{skaras2021}.
With $L$ stages and Gegenbauer parameter $\nu_c$, scaled to a real interval
$[-\beta,0]$, the undamped stage polynomials are
\begin{equation}
\label{eq:chain-stage-undamped}
  R_j^{(c)}(z)=\frac{C_j^{(\nu_c)}(1+v_1z)}{C_j^{(\nu_c)}(1)},\qquad
  v_1=\frac{2}{\beta},\qquad j=0,\ldots,L-1,
\end{equation}
so that $w=1+v_1z$ sends $[-\beta,0]$ onto $[-1,1]$.

The chain can be generalized to include a shift for damping as in RKC methods~\cite{vdh1980,sommeijer1998}, which is not part of the usual RKG construction.
Normalized orthogonal polynomials oscillate on $[-1,1]$, and without damping the
stability polynomial has interior points at which $|R|$ is only marginally below
one.
Shifting the evaluation point from $1$ to
\begin{equation}
\label{eq:chain-damping}
  v_0=1+\frac{\varepsilon_c}{L^{2}},\qquad \varepsilon_c\ge 0,
\end{equation}
pulls those interior maxima strictly below one, at the cost of a slightly shorter
real reach for the same degree.
The damped stage polynomials are then
\begin{equation}
\label{eq:chain-stage}
  R_j^{(c)}(z)=\frac{C_j^{(\nu_c)}(v_0+v_1z)}{C_j^{(\nu_c)}(v_0)},\qquad
  v_1=\frac{1+v_0}{\beta},\qquad j=0,\ldots,L-1,
\end{equation}
which recovers~\eqref{eq:chain-stage-undamped} when $\varepsilon_c=0$.

Alongside the real reach~$\beta$, write $\delta$ for the \emph{largest imaginary
half-height} of a stability polynomial~$R$ over that interval, the largest height
such that every vertical line through $x\in[-\beta,0]$ remains stable at least up to
$\pm\delta$,
\begin{equation}
\label{eq:delta}
  \delta=\max\bigl\{d\ge 0:\ |R(x+iy)|\le 1\text{ for all }x\in[-\beta,0]\text{ and }|y|\le d\bigr\}.
\end{equation}
An \emph{interior damping level} $\eta\in(0,1)$ is a design constraint on how large
$|R|$ may be in the interior of $[-\beta,0]$. Away from the origin, one requires
$|R(x)|\le\eta$ on $[-\beta,-x_0]$ for a cutoff $x_0>0$.
Near the origin, $|R(x)|\approx e^{x}>\eta$, so the bound is imposed only for
$x\le -x_0$ (requiring $x_0>-\log\eta$).
Imposing~$\eta$ keeps interior extrema safely below one, which prevents~$\delta$ from
becoming very small and leaves some stability on the imaginary axis.
The shift~$\varepsilon_c$ and the level~$\eta$ are different controls.
The shift enters the formula through~$v_0$ and changes the shape of~$|R|$.
The level~$\eta$ restricts the admissible
reach~$\beta$ and, through~\eqref{eq:delta}, keeps a useful imaginary half-height,
which is important for mixed hyperbolic--parabolic problems.

Substituting $w=v_0+v_1z$ into~\eqref{eq:gegenbauer} and normalizing gives a
three-term Runge--Kutta recurrence
\begin{equation}
\label{eq:chain-recurrence}
  Y_0=y_n,\qquad
  Y_1=Y_0+\tilde\mu_1\,h f(Y_0),\qquad
  Y_j=\mu_jY_{j-1}+\kappa_jY_{j-2}+\tilde\mu_j\,h f(Y_{j-1}),\quad j=2,\ldots,L-1,
\end{equation}
with coefficients
\begin{equation}
\label{eq:chain-coefficients}
\begin{aligned}
  \tilde\mu_1=\frac{v_1}{v_0},\qquad
  \mu_j=\frac{2(j-1+\nu_c)}{j}\,v_0\,\frac{C_{j-1}^{(\nu_c)}(v_0)}{C_j^{(\nu_c)}(v_0)},\\
  \kappa_j=-\frac{j+2\nu_c-2}{j}\,\frac{C_{j-2}^{(\nu_c)}(v_0)}{C_j^{(\nu_c)}(v_0)},\qquad
  \tilde\mu_j=\frac{2(j-1+\nu_c)}{j}\,v_1\,\frac{C_{j-1}^{(\nu_c)}(v_0)}{C_j^{(\nu_c)}(v_0)},
\end{aligned}
\end{equation}
so that $Y_j=R_j^{(c)}(z)\,y_n$ on the linear test equation and $\mu_j+\kappa_j=1$.
Each stage approximates the solution at the stage time
\begin{equation}
\label{eq:chain-times}
  c_j=R_j^{(c)\prime}(0)=v_1\,\frac{C_j^{(\nu_c)\prime}(v_0)}{C_j^{(\nu_c)}(v_0)}=2\nu_c v_1\,\frac{C_{j-1}^{(\nu_c+1)}(v_0)}{C_j^{(\nu_c)}(v_0)},
\end{equation}
the $j$th row sum of the Butcher matrix of the chain.

The following propositions relate a single chain to Definition~\ref{def:physical-stages}.
Proposition~\ref{prop:chain-stages} shows that the stages are real-valued, consistent, and bounded on the stability interval, which is conditions~(i)--(iii) of the definition.
Proposition~\ref{prop:chain-strong} builds on this by identifying each stage as the full step of a first-order Gegenbauer method and by proving that the stage times are strictly increasing, so that condition~(iv) holds and the chain is physical and monotone in time.

\begin{proposition}
  \label{prop:chain-stages}
  Let $\nu_c>0$ and $\varepsilon_c\ge 0$.
  Then every stage polynomial~\eqref{eq:chain-stage} has real coefficients, satisfies
  $R_j^{(c)}(z)=1+c_jz+O(z^2)$ with $c_j>0$ for $j\ge 1$, and satisfies
  $|R_j^{(c)}(x)|\le 1$ for all $x\in[-\beta,0]$.
  In particular, every stage of a chain is physical and monotone in time in the sense
  of Definition~\ref{def:physical-stages}.
\end{proposition}

\begin{proof}
The coefficients are real by construction, and $c_j>0$ since $v_0\ge 1$ lies at or
beyond the largest zero of $C_j^{(\nu_c)}$, where both $C_j^{(\nu_c)}$ and its
derivative are positive.
For $x\in[-\beta,0]$, the argument $v_0+v_1x$ lies in $[-1,v_0]$.
On $[-1,1]$, $|C_j^{(\nu_c)}(w)|\le C_j^{(\nu_c)}(1)$ for $\nu_c>0$, and on
$[1,v_0]$, $C_j^{(\nu_c)}$ is positive and increasing, hence
$|C_j^{(\nu_c)}(w)|\le C_j^{(\nu_c)}(v_0)$ on $[-1,v_0]$.
For condition~(iv), the stage $Y_j$ is computed from $Y_{j-1}$ and $Y_{j-2}$,
and~\eqref{eq:chain-recurrence} gives $c_j=\mu_jc_{j-1}+\kappa_jc_{j-2}+\tilde\mu_j$.
Since $\mu_j+\kappa_j=1$,
\[
  c_j=c_{j-1}+\kappa_j\,(c_{j-2}-c_{j-1})+\tilde\mu_j.
\]
For $j\ge 2$, $\kappa_j\le 0$ and $\tilde\mu_j>0$ by~\eqref{eq:chain-coefficients}, so
$c_{j-1}>c_{j-2}$ implies $c_j>c_{j-1}>c_{j-2}$, and $c_1=\tilde\mu_1>0=c_0$ starts
the induction.
\end{proof}

\begin{proposition}
  \label{prop:chain-strong}
  Let $\nu_c>0$ and $\varepsilon_c\ge 0$, and let $R_j^{(c)}$ and $c_j$ be the stage
  polynomials and stage times~\eqref{eq:chain-stage} and~\eqref{eq:chain-times} of a
  chain.
  Then:
  \begin{enumerate}
  \item[(a)] for each $j$, truncating the chain after $j$ stages gives a
  (possibly damped) first-order Gegenbauer method whose step is $Y_j$, and
  $|R_j^{(c)}(x)|\le 1$ for all $x\in[-\beta,0]$; if $\varepsilon_c>0$, moreover
  $|R_j^{(c)}(x)|\le C_j^{(\nu_c)}(1)/C_j^{(\nu_c)}(v_0)<1$ for
  $x\in[-\beta,\,-(v_0-1)/v_1]$;
  \item[(b)] the stage times are strictly increasing, $0=c_0<c_1<\cdots<c_{L-1}$; and
  \item[(c)] the stage times satisfy $c_j\le v_1\,j(j+2\nu_c)/(2\nu_c+1)$; in
  particular, $c_j\le 1$ for all $j\le L-1$ whenever
  $\beta\ge(1+v_0)(L-1)(L-1+2\nu_c)/(2\nu_c+1)$.
  \end{enumerate}
\end{proposition}

\begin{proof}
The argument uses standard properties of the Gegenbauer polynomials~\cite{szego1975,dlmf}.
(a) The truncated recurrence has stability function $R_j^{(c)}$ by construction.
For $x\in[-\beta,0]$, $v_0+v_1x\in[-1,v_0]$, and the bound follows as in
Proposition~\ref{prop:chain-stages}.
If $x\le-(v_0-1)/v_1$, then $v_0+v_1x\in[-1,1]$, where
$|C_j^{(\nu_c)}(w)|\le C_j^{(\nu_c)}(1)<C_j^{(\nu_c)}(v_0)$ for $v_0>1$.

(b) Let $x_1^{(j)}>\cdots>x_j^{(j)}$ be the zeros of $C_j^{(\nu_c)}$, all in
$(-1,1)$.
Then $C_j^{(\nu_c)\prime}(v_0)/C_j^{(\nu_c)}(v_0)=\sum_{i=1}^{j}1/(v_0-x_i^{(j)})$,
with every term positive since $v_0\ge 1$.
By the interlacing of the zeros of orthogonal polynomials,
$x_i^{(j+1)}>x_i^{(j)}$ for $i=1,\ldots,j$, so each of the first $j$ terms increases
from $j$ to $j+1$, and an additional positive term is added; hence $c_{j+1}>c_j$.

(c) Since $v_0\ge 1$, each term satisfies $1/(v_0-x_i^{(j)})\le 1/(1-x_i^{(j)})$, so
$c_j\le v_1\,C_j^{(\nu_c)\prime}(1)/C_j^{(\nu_c)}(1)=v_1\,j(j+2\nu_c)/(2\nu_c+1)$,
using $C_j^{(\nu_c)}(1)=(2\nu_c)_j/j!$ and~\eqref{eq:gegenbauer-derivative}.
The last statement follows from $v_1=(1+v_0)/\beta$ and the monotonicity in~(b).
\end{proof}

Proposition~\ref{prop:chain-strong} strengthens Definition~\ref{def:physical-stages}
for each chain: every stage is the full step of a classical first-order stabilized
method (with optional RKC-style damping), and the stages advance monotonically in
time.

\subsection{The target stability polynomial}
\label{sec:crkg-target}

To construct methods of higher order, several chains must be combined, so that the
free weights in~\eqref{eq:crkg-update}--\eqref{eq:crkg-stability} provide enough
freedom to satisfy the higher-order conditions.
Before constructing these chains, the stability polynomial~\eqref{eq:crkg-stability}
of the combined method should first be determined, since it determines the design of
the chains.
As written in~\eqref{eq:crkg-stability}, that polynomial is the result of several
different chains; the polynomial it is required to match is called the \emph{target}
stability polynomial.

The target is prescribed to be a Gegenbauer stability polynomial built from the same
ingredients as the single-chain stage polynomials~\eqref{eq:chain-stage}, with
Gegenbauer parameter~$\nu$, real reach~$\beta$, and RKC-style damping shift~$\varepsilon$,
together with an order~$p$.
Following O'Sullivan~\cite{osullivan2019}, the target of order~$p$ is assembled from
$p$ Gegenbauer components whose degrees are successive multiples of a positive integer
$M$, degrees $M,2M,\ldots,pM$.
$M$ is called the \emph{segment count} because each component advances the degree by one
segment of length~$M$. Increasing $M$ lengthens the real reach while the number
of components stays equal to the order~$p$.
With $L=pM$, the target is
\begin{equation}
\label{eq:target}
  R_{p,M}(z)=d_0+2\sum_{k=1}^{p}d_k\,C_{kM}^{(\nu)}\bigl(w(z)\bigr),\qquad
  w(z)=w_0+\frac{1+w_0}{\beta}\,z,\qquad
  w_0=1+\frac{\varepsilon}{M^{2}},
\end{equation}
with coefficients fixed by $R_{p,M}(z)=e^{z}+O(z^{p+1})$.
The shift~$\varepsilon$ and interior damping level~$\eta$ are the two damping notions
of Section~\ref{sec:crkg-chain}. $\beta$ is taken as large as possible subject to
$|R_{p,M}(x)|\le 1$ on $[-\beta,0]$ and the interior bound with~$\eta$.
This target polynomial will be the stability function of the resulting CRKG method
once the multiple chains and weights are chosen within their constraints.

The default parameters used throughout this paper are $\nu=1/2$, $\varepsilon=1/2$,
$\eta=0.8$, and $x_0=1/2$. These give
$\beta\approx 0.30\,L^{2}$ for $p=3$ and $\beta\approx 0.21\,L^{2}$ for $p=4$.
The cutoff need only satisfy $x_0>-\log\eta$. Instead of setting $x_0=-\log\eta$, the
fixed value $x_0=1/2$ is taken, which for $\eta=0.8$ lies comfortably above
$-\log 0.8\approx 0.223$.
This choice of target parameters is somewhat heuristic, but we have found it to be
the best for our goals.
With $\nu=1/2$, the Gegenbauer polynomials in~\eqref{eq:target}
return the classical Legendre polynomials.

\subsection{Shared stages}
\label{sec:crkg-sharing}

If $n_c$ chains are evaluated independently, each of length~$L$, the step costs
$n_c L-(n_c-1)$ right-hand side evaluations, since only $f(y_n)$ is common, resulting in nearly
$n_c$ times the cost of a single chain.
To reduce this cost, the chains may share their first $k$ stages past the first stage.
One chain, the \emph{trunk}, is first evaluated in full.
The remaining chains, the \emph{branches}, reuse the first $k+1$ trunk stages
$Y_{1,0},\ldots,Y_{1,k}$ and their right-hand sides, then continue with their own
recurrence coefficients,
\begin{equation}
\label{eq:branch-restart}
  Y_{c,k+1}=\mu^{(c)}_{k+1}Y_{1,k}+\kappa^{(c)}_{k+1}Y_{1,k-1}+\tilde\mu^{(c)}_{k+1}\,h f(Y_{1,k}),
\end{equation}
with $Y_{c,1}=y_n+\tilde\mu^{(c)}_1\,h f(y_n)$ when $k=0$, followed
by~\eqref{eq:chain-recurrence} for the remaining stages.
A design with $k$ shared stages costs
\begin{equation}
\label{eq:rhs-per-step}
  L+(n_c-1)(L-1-k)
\end{equation}
right-hand side evaluations per step, so each shared stage saves one evaluation per
branch.
Sharing reduces the number of free stages (and thus the degrees of freedom in the
weights), but it is the mechanism that keeps multi-chain methods practical.
A perturbation of a shared stage reaches every branch, so larger~$k$ tends to increase
internal amplification; how large a~$k$ to take is fixed when the chains are designed
below.

\subsection{Design of the chains}
\label{sec:crkg-chains}

The $n_c$ chains and their weights are now designed so that the method meets the goals
set out above: physical, monotone-in-time stages, the target
polynomial~\eqref{eq:target} with the default parameters, low internal amplification, and as large a shared-stage
count~$k$ as is compatible with a reasonable amplification.
Every stage of every chain already satisfies Definition~\ref{def:physical-stages} by
Propositions~\ref{prop:chain-stages} and~\ref{prop:chain-strong}. What remains is to
choose the chain parameters and the weights in~\eqref{eq:crkg-update} so
that~\eqref{eq:crkg-stability} matches $R_{p,M}$ and the higher-order conditions
hold, at low $\mathcal{A}$ and high~$k$.

In addition to $(\nu_c,\varepsilon_c)$, each chain may introduce a \emph{coverage}
$\sigma_c\ge 1$.
Coverage replaces $v_1=(1+v_0)/\beta$ in~\eqref{eq:chain-stage} by
\begin{equation}
\label{eq:coverage}
  v_1=\frac{1+v_0}{\sigma_c\,\beta},
\end{equation}
so that the chain is scaled to the longer interval $[-\sigma_c\beta,0]$.
On the design interval $[-\beta,0]$ the stages then stay deeper inside the region
where the Gegenbauer polynomials are bounded. Larger coverage damps the stages more
strongly and lowers amplification, at the expense of shorter stage times
$c_j\propto 1/\sigma_c$.
Different chains may use different coverages, and the resulting stability polynomial
is still the target~\eqref{eq:target}, with its reach $\beta$ unchanged by $\sigma_c$.
With $\sigma_c\ge 1$, the bounds of Propositions~\ref{prop:chain-stages}
and~\ref{prop:chain-strong} continue to hold on $[-\beta,0]$, with the natural
interval in part~(a) of Proposition~\ref{prop:chain-strong} extended to
$[-\sigma_c\beta,0]$.
The role of greater coverage is to reduce amplification without sacrificing accuracy
relative to the base variants later in the design of CRKG.

\paragraph{Order conditions}
Write a Runge--Kutta method in Butcher form with matrix $A$, weights $\mathbf{b}$,
and stage times $\mathbf{c}=A\mathbf{1}$.
For each rooted tree~$t$ with $|t|$ vertices, the method has an elementary weight
$\mathbf{b}^{T}\Phi(t)$, and order~$p$ requires $\mathbf{b}^{T}\Phi(t)=1/\gamma(t)$ for
every tree with $|t|\le p$~\cite[Sec.~II.2]{hairer1993},~\cite{butcher2016}.
The conditions up to order four are listed in Table~\ref{tab:order-conditions}.

\begin{table}[H]
\centering
\caption{Order conditions up to order four~\cite{butcher2016}.}
\label{tab:order-conditions}
\begin{tabular}{cllcc}
\toprule
Order & Tree & Condition & $1/\gamma(t)$ & Tall \\
\midrule
1 & $\tau$ & $\mathbf{b}^{T}\mathbf{1}$ & $1$ & yes \\
2 & $[\tau]$ & $\mathbf{b}^{T}\mathbf{c}$ & $1/2$ & yes \\
3 & $[\tau,\tau]$ & $\mathbf{b}^{T}\mathbf{c}^{2}$ & $1/3$ & no \\
3 & $[[\tau]]$ & $\mathbf{b}^{T}A\mathbf{c}$ & $1/6$ & yes \\
4 & $[\tau,\tau,\tau]$ & $\mathbf{b}^{T}\mathbf{c}^{3}$ & $1/4$ & no \\
4 & $[\tau,[\tau]]$ & $\mathbf{b}^{T}(\mathbf{c}\,A\mathbf{c})$ & $1/8$ & no \\
4 & $[[\tau,\tau]]$ & $\mathbf{b}^{T}A\mathbf{c}^{2}$ & $1/12$ & no \\
4 & $[[[\tau]]]$ & $\mathbf{b}^{T}A^{2}\mathbf{c}$ & $1/24$ & yes \\
\bottomrule
\end{tabular}
\end{table}

The \emph{tall} trees have elementary weights equal to the coefficients of $z^{k}$ in
the stability function~\cite{butcher2016}.
If $R(z)=e^{z}+O(z^{p+1})$, all tall-tree conditions up to order~$p$ hold
automatically.
Every condition up to order two is tall, so matching the target alone yields
nonlinear order two.
The remaining conditions are non-tall and are what the extra freedom of the chain
design must satisfy.
Order three needs one non-tall condition, and order four needs four.

\paragraph{Weights}
For fixed chains,~\eqref{eq:crkg-stability} is linear in the weights.
Since both $R$ and $R_{p,M}$ have degree~$L$ and equal one at $z=0$, they coincide if
they agree at $L$ distinct nonzero points.
The matching $R(x_i)=R_{p,M}(x_i)$ is imposed at the Chebyshev
points~\cite{trefethen2019}
$x_i=-\tfrac{\beta}{2}\bigl(1-\cos\tfrac{(2i-1)\pi}{2L}\bigr)$, $i=1,\ldots,L$,
so that the resulting stability polynomial~\eqref{eq:crkg-stability} matches the target~\eqref{eq:target}.
The non-tall order conditions are likewise linear in the weights.
Collecting the weights in a vector $\mathbf{b}$,
\begin{equation}
\label{eq:weight-system}
  \begin{bmatrix} S \\ T \end{bmatrix}\mathbf{b}=\begin{bmatrix} \bigl(R_{p,M}(x_i)-1\bigr)_{i=1}^{L} \\ \bigl(1/\gamma(t)\bigr)_{t\in\mathcal{T}_p} \end{bmatrix},
\end{equation}
where the rows of $S$ evaluate~\eqref{eq:crkg-stability} at the $x_i$, and the rows of
$T$ are the elementary weights of the non-tall trees $\mathcal{T}_p$ with $|t|\le p$.
With a single chain the $L$ stability conditions fix the $L$ weights uniquely.
Each additional chain adds $L-1$ degrees of freedom (its first stage coincides with
$y_n$), which are used for the non-tall conditions.
Order three uses $n_c=2$ chains and order four uses $n_c=3$. For order four, two
chains already provide enough degrees of freedom, but a third allows substantially
lower internal amplification to meet the design goals.

The system~\eqref{eq:weight-system} is underdetermined, and the remaining freedom is
used to keep the weights and the leading error small.
Let $E$ have the elementary weights of the trees with $|t|=p+1$ as rows and
$\mathbf{e}=(1/\gamma(t))_{|t|=p+1}$.
Among solutions of~\eqref{eq:weight-system}, the objective
\begin{equation}
\label{eq:weight-objective}
  \|\mathbf{b}\|^{2}+\|E\mathbf{b}-\mathbf{e}\|^{2}
\end{equation}
is minimized via a truncated SVD of~\eqref{eq:weight-system}~\cite{golub2013,bjorck1996},
followed by one residual correction of the particular solution.
A design is accepted if the relative residual is below $10^{-10}$. Write
$\mathcal{E}=\|E\mathbf{b}-\mathbf{e}\|$ for the resulting leading error.

\paragraph{Selecting the chains and the shared-stage count}
The chains are drawn from a pool $\mathcal{V}$ of variants with
$\nu_c\in\{1/2,1,2,3\}$, $\varepsilon_c\in\{1/2,2,10,30\}$, and $\sigma_c\in\{1,3/2,2\}$
($48$ variants).
Every combination of $n_c$ variants from the pool is tried, the weights are solved and $\mathcal{A}$ and
$\mathcal{E}$ are recorded.
An accuracy bar $\bar{\mathcal{E}}$ is taken from a smaller base pool
$\mathcal{V}_0\subset\mathcal{V}$ with $\nu_c\in\{1/2,1,2\}$, the same $\varepsilon_c$ values as $\mathcal{V}$, and $\sigma_c=1$; among those
combinations, the one that minimizes $\mathcal{A}\,\mathcal{E}$ sets $\bar{\mathcal{E}}$.
In this subset, $\nu_c=1/2$ and $\nu_c=1$ recover the Legendre and Chebyshev polynomials of the second kind, with $\nu_c=2$ the next Gegenbauer parameter, while $\sigma_c=1$ scales each chain exactly to the target interval~$[-\beta,0]$.
Extra coverage and the larger $\nu_c=3$ are excluded from $\mathcal{V}_0$ because both damp the stages more strongly and shorten stage times, which tends to raise the leading error, whereas all damping shifts $\varepsilon_c$ are retained since damping can improve accuracy.
The lowest error of a chain from this subset is used as an error benchmark.
All combinations with $\mathcal{E}\le\bar{\mathcal{E}}$ are kept as candidates,
done so that the remaining candidates are within a reasonable range of accuracy, and
$\mathcal{A}_{\min}$ denotes the smallest amplification among these candidates.

The shared-stage count is then chosen among these candidates.
For an amplification factor~$\alpha\ge 1$, the largest~$k$ is taken for which some
candidate, with some chain as trunk, still meets the order conditions,
$\mathcal{E}\le\bar{\mathcal{E}}$, and $\mathcal{A}\le\alpha\,\mathcal{A}_{\min}$.
At that~$k$, the design with lowest $\mathcal{A}$ is kept.
With $\alpha=1$ this is the lowest-amplification design; larger $\alpha$ allows more
sharing at the cost of higher amplification.
We fix $\alpha=5$, which keeps amplification in a reasonable range while reducing the
cost per step. We find that above this default, further sharing brings diminishing
returns at the cost of significantly higher amplification.

\paragraph{Refinement}
Once the trunk, branches, and~$k$ are fixed, the coefficients are refined with the following procedure.
With the chains held fixed, the weights are recomputed by a linear program that
minimizes internal amplification subject to~\eqref{eq:weight-system} and to the
requirement that no next-order error component exceeds its value in the discrete
design~\cite{boyd2004}:
\begin{equation}
\label{eq:lp-weights}
  \min_{\mathbf{b},\,\tau}\; \tau
  \quad\text{subject to}\quad
  \eqref{eq:weight-system},\quad
  |\Theta\mathbf{b}|\le \tau\,\mathbf{1},\quad
  |E\mathbf{b}-\mathbf{e}|\le\boldsymbol{\delta}_{\mathrm{e}},
\end{equation}
where the rows of $\Theta$ are the responses of the step output to unit stage
perturbations on $[-\beta,0]$, and $\boldsymbol{\delta}_{\mathrm{e}}$ holds the
componentwise next-order residuals of the discrete design.
The optimization problem~\eqref{eq:lp-weights} is solved with HiGHS~\cite{huangfu2018}.
The parameters $(\nu_c,\varepsilon_c,\sigma_c)$ are then adjusted by a short
Nelder--Mead search~\cite{nelder1965} in $\bigl(\log\nu_c,\,\log\varepsilon_c,\,s_c\bigr)$
with $\sigma_c=1+e^{s_c}$, evaluating~\eqref{eq:lp-weights} at each trial point, and the
result is kept if $\mathcal{A}$ decreases.
The combinatorial choice of trunk, branches, and~$k$, and the target~$R_{p,M}$, are
left unchanged; only the continuous chain parameters and the weights~$\mathbf{b}$ may
change in this refinement step.

\begin{remark}\label{rem:embedded-weights}
Each method produces embedded weights
$\hat{\mathbf{b}}$ of order $p-1$ on the same stages.
The top coefficient of the target is perturbed to $\hat d_p=d_p(1+\kappa_{\mathrm{emb}})$,
and the remaining coefficients are fixed by the order-$p-1$ conditions
\begin{equation}
\label{eq:embedded-order}
  \sum_{k=1}^{p-1}\hat d_k\,\frac{d^{n}C_{kM}^{(\nu)}}{dw^{n}}(w_0)=\frac{1}{2}\Bigl(\frac{\beta}{1+w_0}\Bigr)^{n}-\hat d_p\,\frac{d^{n}C_{pM}^{(\nu)}}{dw^{n}}(w_0),\qquad n=1,\ldots,p-1,
\end{equation}
with $\hat d_0$ from $\hat R_{p,M}(0)=1$.
This provides support for error estimation and adaptive step-size control.
\end{remark}

\subsection{Summary of the construction}
\label{sec:crkg-algorithm}

The objects and default parameters of the construction are summarized in Table~\ref{tab:crkg-defaults}.
The flow of the construction is summarized in Algorithm~\ref{alg:crkg-design}.
The construction for each $M$ is computed once and stored for all
future use.

\begin{table}[H]
\centering
\caption{Objects and default parameters of the CRKG construction.}
\label{tab:crkg-defaults}
\begin{tabular}{@{}l>{\raggedright\arraybackslash}p{0.62\linewidth}@{}}
\toprule
Item & Definition / default \\
\midrule
\multicolumn{2}{@{}l}{\emph{Objects}} \\
Target $R_{p,M}$ & prescribed stability function of the full step (reach $\beta$, linear order $p$) \\
Chain & a single Gegenbauer three-term recurrence (physical, monotone-in-time stages) \\
Chain combination & several chains with a weighted update (nonlinear order $p$) \\
Trunk and branches & chains that share their first stages (reduced cost per step) \\
\midrule
\multicolumn{2}{@{}l}{\emph{Defaults}} \\
$\nu$ & Gegenbauer parameter of the target: $1/2$ \\
$\varepsilon$ & RKC-style damping shift of the target: $1/2$ \\
$\eta$ & interior damping level of the target: $0.8$ \\
$x_0$ & cutoff for the interior bound: $1/2$ \\
$n_c$ & number of chains: $2$ ($p=3$), $3$ ($p=4$) \\
$\mathcal{V}$ & pool of chain variants $(\nu_c,\varepsilon_c,\sigma_c)$:\newline
  $\nu_c\in\{1/2,1,2,3\}$,\newline
  $\varepsilon_c\in\{1/2,2,10,30\}$,\newline
  $\sigma_c\in\{1,3/2,2\}$ \\
$\alpha$ & amplification factor for shared-stage selection: $5$ \\
\bottomrule
\end{tabular}
\end{table}

\begin{algorithm}[H]
\caption{CRKG design for given order~$p$ and segment count~$M$.}
\label{alg:crkg-design}
\begin{algorithmic}[1]
\State Build the target~$R_{p,M}$ with the defaults of Table~\ref{tab:crkg-defaults}; take~$\beta$ as large as possible subject to~$|R|\le 1$ and the interior bound with~$\eta$
\State For each $n_c$-subset of~$\mathcal{V}$: solve~\eqref{eq:weight-system}--\eqref{eq:weight-objective}; record~$\mathcal{A}$ and~$\mathcal{E}$
\State Set the accuracy bar~$\bar{\mathcal{E}}$ from the base variants~$\mathcal{V}_0$; keep candidates with~$\mathcal{E}\le\bar{\mathcal{E}}$
\State Take the largest shared-stage count~$k$ with~$\mathcal{A}\le\alpha\,\mathcal{A}_{\min}$; among those, keep the design of lowest~$\mathcal{A}$
\State Refine the weights by~\eqref{eq:lp-weights} and the continuous chain parameters by Nelder--Mead
\State Compute the embedded weights~$\hat{\mathbf{b}}$ of order~$p-1$
\end{algorithmic}
\end{algorithm}

\begin{remark}
\label{rem:construction-params}
Selecting the chains, shared stages, and weights is a constrained optimization aimed
at low right-hand side cost and low internal amplification subject to the order and
stability conditions.
A discrete search followed by continuous refinement is used. Other formulations, such
as convex optimization for stability polynomials~\cite{ketcheson2012}, may yield still
better designs.
\end{remark}

\section{Step evaluation}
\label{sec:computation}

Once a design has been constructed as in Section~\ref{sec:crkg}, it is stored as a short list of coefficients: the recurrence coefficients $(\mu_j,\kappa_j,\tilde\mu_j)$ of each chain, the number of shared stages~$k$, the weights $b_{c,j}$ and $\hat b_{c,j}$, and the stage times $c_{c,j}$.
Evaluating a step is then a standard orthogonal-polynomial recurrence, of the same form as the classical recurrence-based methods of Section~\ref{sec:prelim-orthogonal}, applied once along the trunk and once along each branch.

A step is evaluated in two passes, summarized in Algorithm~\ref{alg:step}.
In the first pass, the trunk is evaluated in full by~\eqref{eq:chain-recurrence}, and at the shared stage~$k$ it retains the states $Y_{1,k}$ and $Y_{1,k-1}$ and the right-hand side $f(Y_{1,k})$ from which the branches restart.
In the second pass, each branch restarts from these by~\eqref{eq:branch-restart} and continues with its own recurrence coefficients.
Since the update~\eqref{eq:crkg-update} is a weighted sum of right-hand side evaluations, each evaluation can be added to the update as soon as it is computed, and the stages themselves need not be kept in memory.
Together with the three-term recurrence, which needs only the two most recent stages of a chain, the storage of a step is therefore independent of the number of stages~$L$.
If a stagewise procedure such as a limiter is used, it is applied to each stage as soon as it is formed, at its stage time $t_n+c_{c,j}h$, and to $y_{n+1}$ at $t_n+h$.

\begin{algorithm}[H]
\caption{A step of a CRKG method from $y_n$ to $y_{n+1}$.}
\label{alg:step}
\begin{algorithmic}[1]
\Require $y_n$, step size $h$, a CRKG design with $n_c$ chains and $k$ shared stages, an optional stagewise procedure $\Lambda$
\Ensure $y_{n+1}$ and, optionally, the embedded error estimate
\Statex \textit{Trunk}
\State $Y_{1,0}\gets y_n$; $Y_{1,1}\gets y_n+\tilde\mu_1h f(y_n)$; $\Lambda(Y_{1,1},t_n+c_{1,1}h)$
\For{$j=1,\ldots,L-2$}
  \State $Y_{1,j+1}\gets\mu_{j+1}Y_{1,j}+\kappa_{j+1}Y_{1,j-1}+\tilde\mu_{j+1}h f(Y_{1,j})$; $\Lambda(Y_{1,j+1},t_n+c_{1,j+1}h)$
\EndFor
\Statex \textit{Branches}
\For{$c=2,\ldots,n_c$}
  \State form $Y_{c,k+1}$ from $Y_{1,k}$, $Y_{1,k-1}$, and $f(Y_{1,k})$ by~\eqref{eq:branch-restart}; $\Lambda(Y_{c,k+1},t_n+c_{c,k+1}h)$
  \For{$j=k+1,\ldots,L-2$}
    \State $Y_{c,j+1}\gets\mu^{(c)}_{j+1}Y_{c,j}+\kappa^{(c)}_{j+1}Y_{c,j-1}+\tilde\mu^{(c)}_{j+1}h f(Y_{c,j})$; $\Lambda(Y_{c,j+1},t_n+c_{c,j+1}h)$
  \EndFor
\EndFor
\Statex \textit{Update}
\State $y_{n+1}\gets y_n+h\sum_{c}\sum_{j}b_{c,j}f(Y_{c,j})$, accumulated as each $f(Y_{c,j})$ is computed; $\Lambda(y_{n+1},t_n+h)$
\State error estimate $h\sum_{c}\sum_{j}(b_{c,j}-\hat b_{c,j})f(Y_{c,j})$, accumulated in the same way
\end{algorithmic}
\end{algorithm}

Each right-hand side $f(Y_{c,j})$ is evaluated once and used both in the recurrence and in the update, so the cost of a step is~\eqref{eq:rhs-per-step} evaluations.

\section{Comparison of CRKG methods with DUMKA and ROCK}
\label{sec:discussion}

In this section we verify that the constructed CRKG methods satisfy the two properties of Section~\ref{sec:prelim}, physical, monotone-in-time stages and low internal amplification, and compare them with the existing methods of the same order: CRKG3 with the third-order DUMKA3 method~\cite{medovikov1998}, and CRKG4 with the fourth-order ROCK4 method~\cite{abdulle2002}.
For a prescribed reach~$\beta$, each method uses its smallest stage count whose real stability interval covers $[-\beta,0]$, and all quantities are evaluated on the linear test equation over $[-\beta,0]$.
For both ROCK4 and DUMKA3 we use the implementation in OrdinaryDiffEq.jl~\cite{rackauckas2017}.

\subsection{Physical, monotone-in-time stages}
\label{sec:discussion-physical}

Table~\ref{tab:physical} reports, for each method and reach, the largest value of $|R_j(x)|$ over all stages and all $x\in[-\beta,0]$, which verifies condition~(iii) of Definition~\ref{def:physical-stages}.
For CRKG, the maximum is taken over all stages of the method, including the stages of the branches, whose boundedness is not covered by Proposition~\ref{prop:chain-stages}.
Every stage of CRKG3 and CRKG4 is bounded by one, so that the branch stages are also physical, and the same holds for ROCK4, whose stages are all real, consistent, and bounded.
DUMKA3, in contrast, has intermediate states that exceed one already at $\beta=50$, and that grow to more than $500$ at $\beta=500$.

\begin{table}[H]
\centering
\caption{Largest stage value $\max_j\max_{x\in[-\beta,0]}|R_j(x)|$ (condition~(iii)).}
\label{tab:physical}
\begin{tabular}{lrrrrc}
\toprule
Method & $\beta$ & $L$ or $s$ & Stage & $\max|R_j|$ & Condition~(iii) \\
\midrule
CRKG3  & 50  & 15 & 4  & $1.0000$ & yes \\
CRKG3  & 200 & 30 & 9  & $1.0000$ & yes \\
CRKG3  & 500 & 42 & 64 & $1.0000$ & yes \\
DUMKA3 & 50  & 15 & 14 & $1.4249$ & no \\
DUMKA3 & 200 & 27 & 26 & $4.6645$ & no \\
DUMKA3 & 500 & 48 & 35 & $5.7117\times10^{2}$ & no \\
\midrule
CRKG4  & 50  & 16 & 32 & $1.0000$ & yes \\
CRKG4  & 200 & 32 & 79 & $1.0000$ & yes \\
CRKG4  & 500 & 52 & 42 & $1.0000$ & yes \\
ROCK4  & 50  & 13 & 1  & $1.0000$ & yes \\
ROCK4  & 200 & 24 & 1  & $1.0000$ & yes \\
ROCK4  & 500 & 38 & 1  & $1.0000$ & yes \\
\bottomrule
\end{tabular}
\end{table}

Condition~(iv) of Definition~\ref{def:physical-stages} can be shown by the construction of the methods.
For CRKG, it holds for every stage, including the stages of the branches, by Proposition~\ref{prop:chain-stages}, and for ROCK4 and DUMKA3 it fails directly from their published coefficients, since the first stage of the finishing procedure of ROCK4 and the first stage of the last substep of DUMKA3 are forward Euler steps of negative step size (Section~\ref{sec:prelim-physical}).

\subsection{Internal amplification}
\label{sec:discussion-amp}

\begin{figure}[H]
\centering
\includegraphics[width=0.7\linewidth]{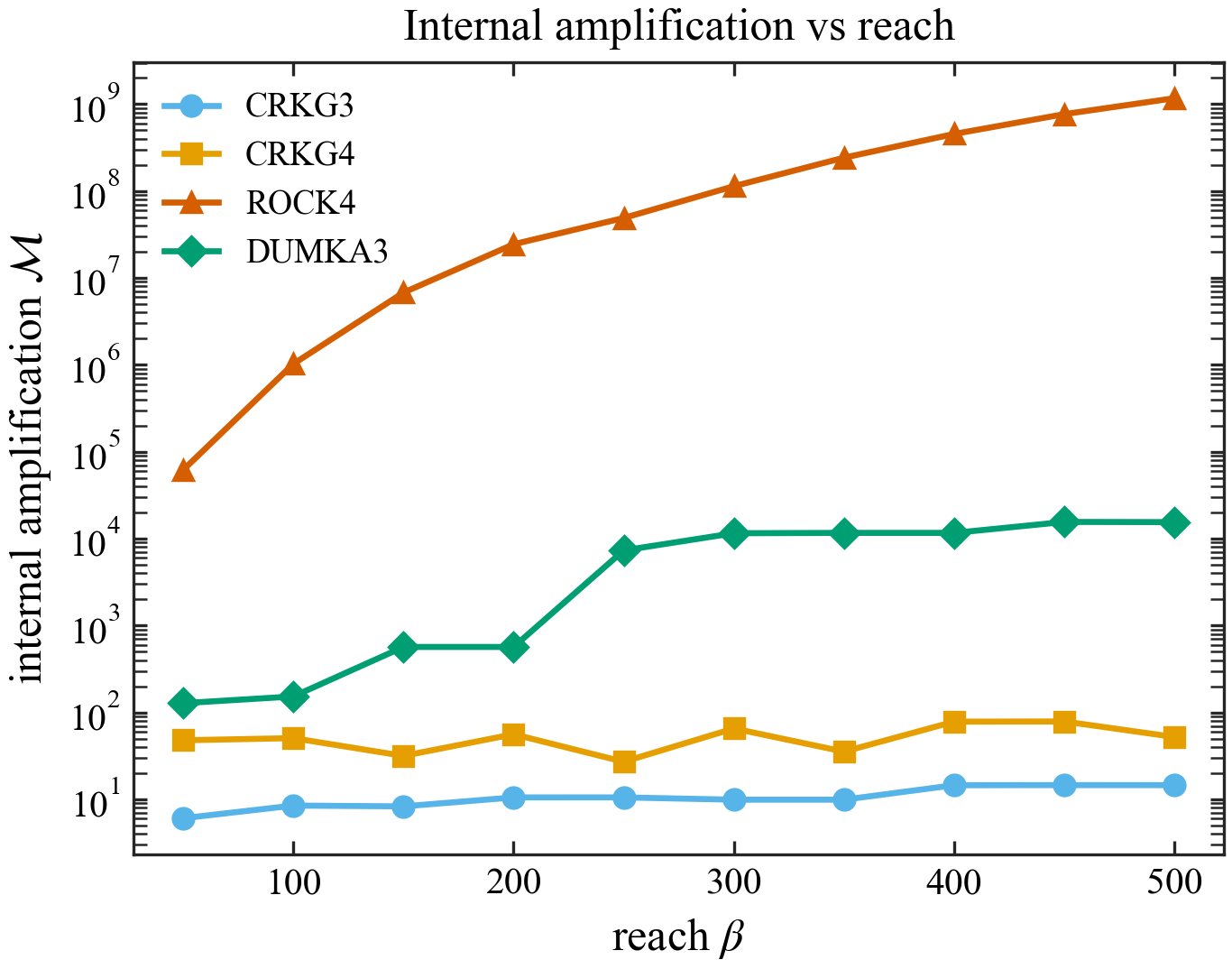}
\caption{Internal amplification $\mathcal{A}$ versus reach $\beta$.}
\label{fig:amp}
\end{figure}

We write $\theta_j(z)$ for the stage amplification factor and $\mathcal{A}$ for the internal amplification of Definition~\ref{def:internal-amp}.
We measure $\mathcal{A}$, following~\cite{ketcheson2014}, for prescribed reaches $\beta$ between $50$ and $500$, increasing the stage count of each method as needed to cover $[-\beta,0]$.
The results are shown in Figure~\ref{fig:amp} and Table~\ref{tab:amp}, together with the number of right-hand side evaluations per step.

\begin{table}[H]
\centering
\caption{Internal amplification $\mathcal{A}$ and RHS evaluations per step $N_f$.}
\label{tab:amp}
\begin{tabular}{rrrrrrrrr}
\toprule
& \multicolumn{2}{c}{CRKG3} & \multicolumn{2}{c}{DUMKA3} & \multicolumn{2}{c}{CRKG4} & \multicolumn{2}{c}{ROCK4} \\
\cmidrule(lr){2-3}\cmidrule(lr){4-5}\cmidrule(lr){6-7}\cmidrule(lr){8-9}
$\beta$ & $\mathcal{A}$ & $N_f$ & $\mathcal{A}$ & $N_f$ & $\mathcal{A}$ & $N_f$ & $\mathcal{A}$ & $N_f$ \\
\midrule
50  & $6.1$ & 24 & $1.3\times10^{2}$ & 15 & $48$ & 28  & $6.2\times10^{4}$ & 13 \\
100 & $8.5$ & 36 & $1.5\times10^{2}$ & 21 & $51$ & 44  & $1.0\times10^{6}$ & 18 \\
200 & $11$ & 49 & $5.7\times10^{2}$ & 27 & $56$ & 62  & $2.4\times10^{7}$ & 24 \\
300 & $9.9$ & 62 & $1.2\times10^{4}$ & 36 & $65$ & 78  & $1.1\times10^{8}$ & 30 \\
400 & $15$ & 68 & $1.2\times10^{4}$ & 36 & $78$ & 94  & $4.5\times10^{8}$ & 34 \\
500 & $15$ & 68 & $1.5\times10^{4}$ & 48 & $52$ & 106 & $1.2\times10^{9}$ & 38 \\
\bottomrule
\end{tabular}
\end{table}

The internal amplification of CRKG3 remains of order $10$ and that of CRKG4 of order $10$--$10^{2}$ over the whole range.
The amplification of DUMKA3 is between one and three orders of magnitude larger than that of CRKG3, and grows from about $10^{2}$ to above $10^{4}$, while that of ROCK4 is between three and seven orders of magnitude larger than that of CRKG4, and grows from about $6\times10^{4}$ to above $10^{9}$.
The amplification of CRKG4 does not vary monotonically with $\beta$, since each reach is covered by a different design, whose chains and number of shared stages are selected independently by the construction of Section~\ref{sec:crkg}.

\begin{remark}[Cost]
\label{rem:cost}
The low amplification of CRKG comes at a cost in right-hand side evaluations as shown in Table~\ref{tab:amp}.
At $\beta=500$, CRKG3 requires $68$ evaluations per step against $48$ for DUMKA3, and CRKG4 requires $106$ against $38$ for ROCK4, so that ROCK4 is roughly a factor of three cheaper per step at the same reach.
This cost stems from the chain combination, which evaluates several chains to obtain the freedom for the nonlinear order conditions, and from the shorter interval of the target polynomial. Stage sharing reduces the former but does not remove it.
This means that, in practice, the cost of a CRKG3 step is around $1.5$--$2\times$ that of a DUMKA3 step, and that of a CRKG4 step is around $2$--$3\times$ that of a ROCK4 step. However, as this section shows, CRKG has important properties that the others do not possess.
In particular, we show later that these properties make CRKG suitable for the stagewise limiting applications that are the focus of this paper, while DUMKA3 and ROCK4 are not.
\end{remark}

Together, Tables~\ref{tab:physical}--\ref{tab:amp} show that CRKG3 and CRKG4 satisfy all of the requirements of Section~\ref{sec:prelim}: their stages are physical and monotone in time, and their internal amplification remains low as the reach increases.
DUMKA3 satisfies neither requirement, and ROCK4, although its stages are physical, steps backward in time and amplifies stage perturbations by factors of $10^{8}$ to $10^{9}$.
CRKG therefore achieves the combination of native high order, physical and monotone-in-time stages, and low internal amplification that the existing higher-order stabilized methods do not, at the price of more right-hand side evaluations per step.

\subsection{Effect of parameters}
\label{sec:discussion-params}

As noted in Remark~\ref{rem:construction-params}, the construction is a constrained optimization over many free parameters, including the target parameters $\nu$, $\varepsilon$, $\eta$, and $x_0$, the chain pool $(\nu_c,\varepsilon_c,\sigma_c)$, the number of chains $n_c$, the amplification factor $\alpha$, and the details of the discrete search and continuous refinement.
We acknowledge that, by itself, this is a significant number of parameters.
The effect of each of these on reach, damping, cost, and internal amplification is already described in the construction of Section~\ref{sec:crkg}; for example, larger $\nu$ or $\varepsilon$ damps the interior of the target at the expense of $\beta$, coverage $\sigma_c$ trades unused reach in a chain for lower amplification, and $\alpha$ trades sharing (hence cost) against amplification.
We do not report a separate parametric study, nor do we detail alternative choices, because we have settled on the defaults and routines of Section~\ref{sec:crkg} as best for our goals---physical, monotone-in-time stages, low internal amplification, and practical efficiency---and we do not recommend changing them.

\section{Numerical results}
\label{sec:numerics}

In this section, numerical results are presented for the CRKG methods.
To begin, convergence is verified on the van der Pol oscillator.
An application of stagewise limiting on a simple one-dimensional advection--diffusion boundary-layer problem then demonstrates the properties of the CRKG methods while showing the failure of DUMKA3 and ROCK4.
The Daru--Tenaud viscous shock tube, the main compressible Navier--Stokes benchmark of this paper, is used next to compare timings of CRKG against adaptive SSPRK and RDPK methods under stagewise limiting, followed by a graded-mesh comparison at $\mathrm{Re}=200$.
Timings use Julia~1.13.4 on a dual-socket AMD {EPYC}~9745 node ($2\times128$ cores, $512$ hardware threads) with $256$ threads pinned to physical cores and NUMA-interleaved memory allocation.

\subsection{Convergence on the van der Pol oscillator}
\label{sec:numerics-vdp}

\begin{figure}[H]
\centering
\includegraphics[width=0.7\linewidth]{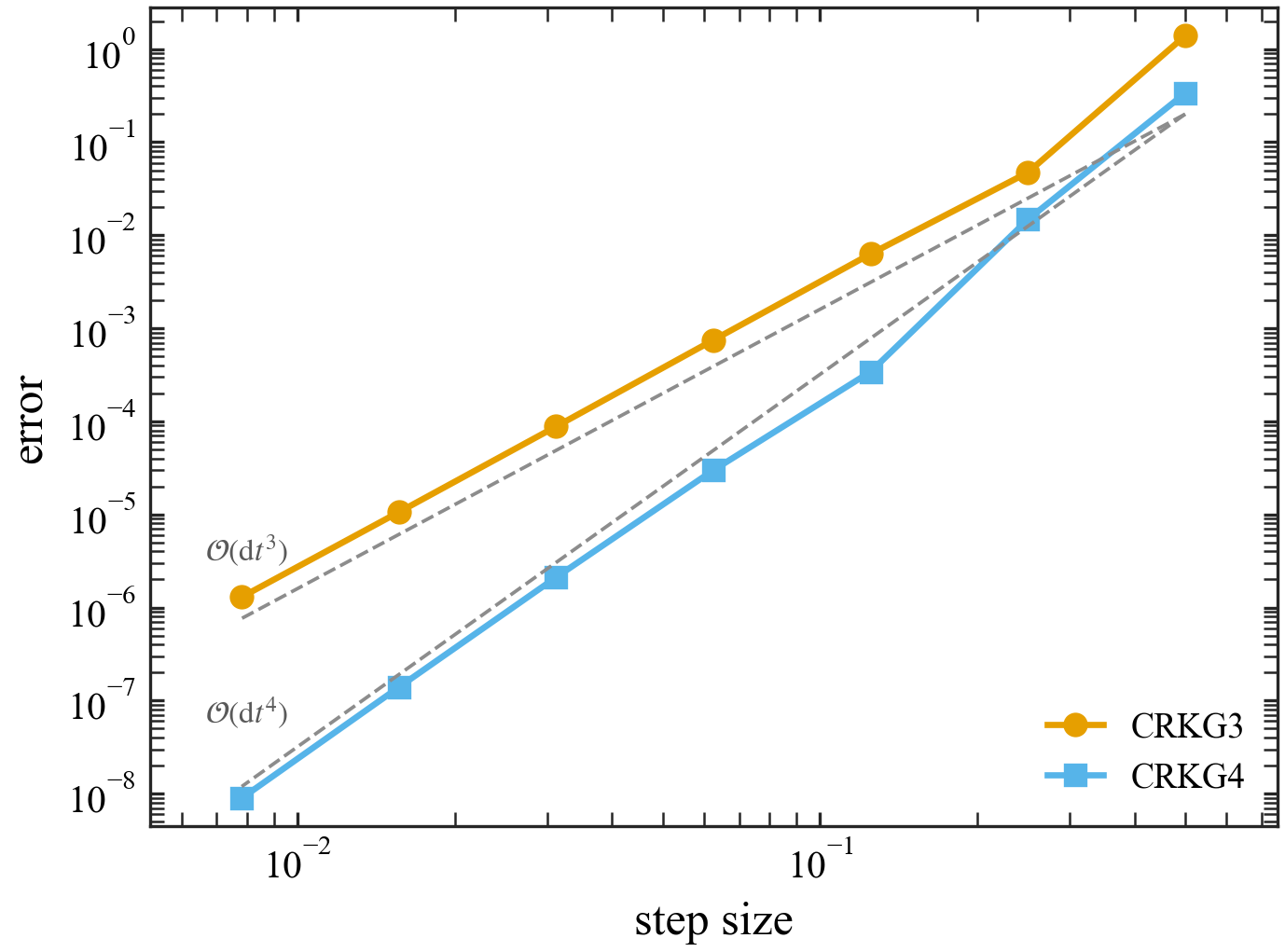}
\caption{Convergence of CRKG3 and CRKG4 on the van der Pol oscillator.}
\label{fig:vdp-conv}
\end{figure}

The van der Pol oscillator~\cite{vanderpol1926,hairer1996}
\begin{equation}
\label{eq:vdp}
  y_1'=y_2,\qquad
  y_2'=\mu(1-y_1^{2})y_2-y_1,
\end{equation}
is a standard nonlinear test problem for stiff and nonstiff ODE solvers.
The stiffness parameter is $\mu=2$, which is mildly stiff, with initial data $(y_1,y_2)=(2,0)$, and the problem is integrated to $t=10$.
A reference solution is computed with the twelfth-order Feagin method~\cite{feagin2012} at a fine step size, and both CRKG3 and CRKG4 are run at fixed segment count $M=8$ over a ladder of step sizes.
Figure~\ref{fig:vdp-conv} shows the absolute error at the final time against the step size.
The observed slopes match the design orders three and four.

\subsection{Application to stagewise limiting}
\label{sec:numerics-ad}

The next application is to stagewise limiting on a simple linear advection--diffusion problem in one dimension,
\begin{equation}
\label{eq:ad1d}
  u_t + u_x = \nu\,u_{xx},
  \qquad x\in[0,1],\quad t\in[0,1],
\end{equation}
with diffusivity $\nu=0.1$.
The initial condition is the constant state $u\equiv 1$, with Dirichlet boundary condition $u=1$ at the inflow $x=0$ and $u=10^{-6}$ at $x=1$.
The problem is discretized with a discontinuous Galerkin spectral element method (DGSEM) of polynomial degree three on a mesh of $16$ elements, using the hyperbolic--parabolic solvers of Trixi.jl~\cite{schlottkelakemper2021trixi}.
At every Runge--Kutta stage the local Zhang--Shu positivity limiter~\cite{zhang2010maximum,zhang2010positivity} is applied together with the global Liu--Zhang optimization-based correction~\cite{liu2024optimization,liu2024simple}, with projection tolerance $10^{-6}$.

\begin{figure}[H]
\centering
\includegraphics[width=\linewidth]{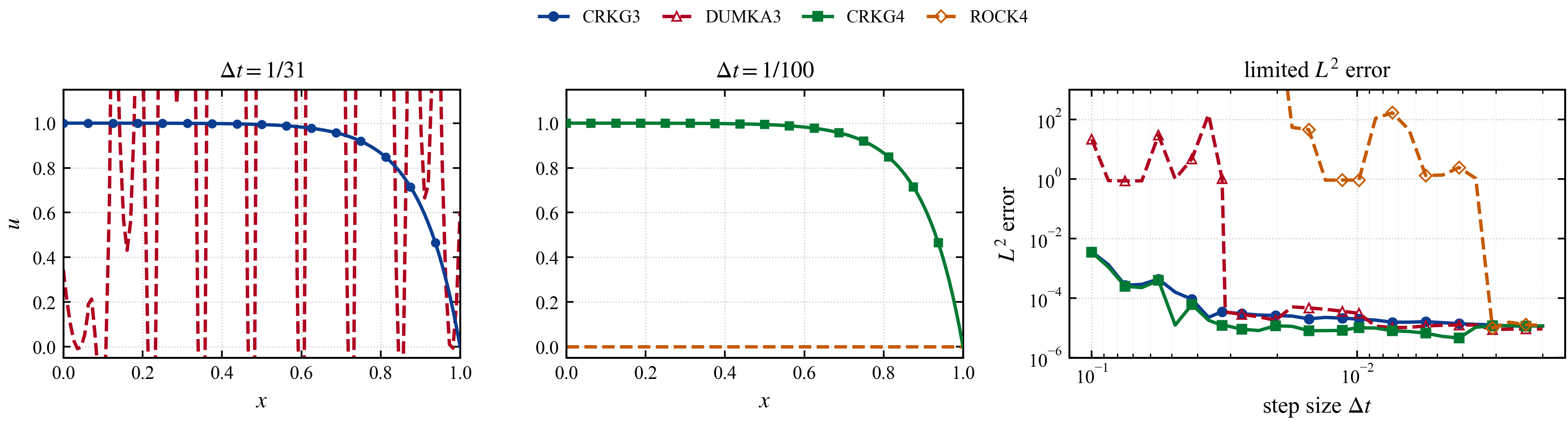}
\caption{Stagewise-limited~\eqref{eq:ad1d} at $t=1$: CRKG3 / DUMKA3 at $\Delta t=1/31$, CRKG4 / ROCK4 at $\Delta t=1/100$, and limited $L^{2}$ error versus $\Delta t$.}
\label{fig:ad1d-limiting}
\end{figure}

Without limiting, both DUMKA3 and ROCK4 are stable on this problem.
With stagewise limiting, both CRKG methods remain physical up to a higher step size, while DUMKA3 and ROCK4 fail until the step is small enough.
Figure~\ref{fig:ad1d-limiting} shows two representative final profiles together with the limited $L^{2}$ error against a fine reference.
At $\Delta t=1/31$, DUMKA3 still oscillates strongly while CRKG3 is physical.
At $\Delta t=1/100$, ROCK4 has been driven to the limiter floor $10^{-6}$ by its amplification and remains a near-constant state, while CRKG4 is physical.
The $L^{2}$ panel summarizes this over the step-size ladder, where CRKG stays near low error for all step sizes, whereas DUMKA3 and ROCK4 remain large until they drop to the same floor.

\subsection{Viscous shock tube of Daru and Tenaud}
\label{sec:numerics-daru}

The main benchmark of this paper is the viscous shock tube of Daru and Tenaud~\cite{daru2001}, a standard test for compressible Navier--Stokes solvers with shocks and wall boundary layers.
The domain is the unit square $[0,1]^{2}$, where $x=1/2$ separates a high-density gas ($\rho_L=120$) from a low-density gas ($\rho_R=1.2$) at rest, with equal temperatures and slip or no-slip walls according to the formulation.
After the diaphragm ruptures, a shock and contact travel into the light gas while a rarefaction expands into the heavy gas. The shock reflects from the end wall and interacts with the contact and with the viscous boundary layers that form along the walls.
The combination of a strong shock, a thin wall boundary layer, and shock--boundary-layer interaction makes the problem stiff in the parabolic modes, sensitive to positivity of density and pressure, and expensive for standard explicit methods whose step size is limited by the viscous CFL condition.
It is therefore taken as the principal demonstration that CRKG retains physical, monotone-in-time stages under stagewise limiting while removing the parabolic step-size restriction.

The compressible Navier--Stokes equations are discretized with a DGSEM of polynomial degree $N=3$ on a Cartesian mesh of $16\,384$ elements, again using Trixi.jl~\cite{schlottkelakemper2021trixi}, at Reynolds number $\mathrm{Re}=200$ and Prandtl number $\mathrm{Pr}=0.73$, and integrated to $t=1$.
At every stage the Zhang--Shu limiter and the Liu--Zhang global correction designed for the compressible Navier--Stokes equations~\cite{zhang2017positivity,liu2024optimization,liu2024simple} are applied.
Four time integrators are compared under this common spatial discretization and limiter.
The third-order strong-stability-preserving method SSPRK43~\cite{spiteri2002,gottlieb2001strong} and the optimized low-storage method RDPK3SpFSAL49~\cite{ranocha2022} are run adaptively, with absolute and relative tolerances $(10^{-6},10^{-4})$ for SSPRK43 and $(10^{-8},10^{-6})$ for RDPK3SpFSAL49.
CRKG3 and CRKG4 select the step size and segment count from CFL estimates of this discretization.
With real and imaginary spectral radii $\rho_r$ and $\rho_i$, real reach $\beta(M)$, and imaginary half-height $\delta(M)$
the admissible step for each~$M$ is
\begin{equation}
\label{eq:cfl-limits}
h(M)=\min\left\{
  \frac{\beta(M)}{\rho_r},\;
  C_i\,\frac{\delta(M)}{\rho_i}
\right\}, 
\end{equation}
and the segment count
\begin{equation}
\label{eq:cfl-choice}
  M^{*}=\operatorname*{arg\,min}_{2\le M\le M_{\max}}\;\frac{N_f(M)}{h(M)},
  \qquad h=h(M^{*}),
\end{equation}
minimizes the cost per unit time, with $N_f(M)$ from~\eqref{eq:rhs-per-step}.
All numerical results in this work use the values $C_i=1.0$ for CRKG3 and $C_i=0.8$ for CRKG4; these represent looser settings which remain physical for this problem.
The spectral-radius estimates used in this work are
\begin{equation}
\label{eq:daru-radii}
\rho_i=\frac{N\sqrt{N+1}}{h_e}\max_{x\in\Omega_h}\sum_{i=1}^{d}\bigl(\lvert v_i\rvert+c_s\bigr), \qquad
\rho_r=\frac{5(N+1)^{4}\,\nu_{\mathrm{eff}}\,d}{4\,h_e^{2}}
+\frac{N(N+1)}{h_e}\max_{x\in\Omega_h}\sum_{i=1}^{d}\bigl(\lvert v_i\rvert+c_s\bigr),
\end{equation}
where $h_e$ is the element width, $c_s=\sqrt{\gamma p/\rho}$ is the local sound speed, and $\nu_{\mathrm{eff}}=\max(\mu/\rho,\,\gamma\mu/(\mathrm{Pr}\,\rho))$ is taken from the minimum density in the current state.
The factor $\tfrac{5}{4}(N+1)^{4}$ and the form of $\phi$ are upper bounds from tests on the spectra of this discretization for $N=1,\ldots,6$.

Table~\ref{tab:daru-timing} reports timings for $N=3$ and, on the same mesh, for $N=4$ as a stress test that stiffens both CFL estimates through~\eqref{eq:daru-radii}.
All runs remain physical through $t=1$.
At $N=3$, CRKG3 is about three times faster than SSPRK43 and CRKG4 is about three times faster than RDPK3SpFSAL49.
At $N=4$, the traditional methods are limited by the early-time boundary-layer parabolic stiffness, while CRKG finishes significantly faster. CRKG3 is about ten times faster than SSPRK43 and CRKG4 about five times faster than RDPK3SpFSAL49.

\begin{table}[H]
\centering
\begin{tabular}{clrrr}
\toprule
$N$ & Method & Steps & RHS evaluations & Wall time [s] \\
\midrule
3 & CRKG3 & $1382$ & $24\,107$ & $57$ \\
3 & CRKG4 & $1152$ & $37\,500$ & $91$ \\
3 & SSPRK43 & $7268$ & $30\,654$ & $165$ \\
3 & RDPK3SpFSAL49 & $6344$ & $63\,559$ & $313$ \\
\midrule
4 & CRKG3 & $2021$ & $41\,714$ & $99$ \\
4 & CRKG4 & $1674$ & $61\,248$ & $143$ \\
4 & SSPRK43 & $28\,083$ & $134\,134$ & $964$ \\
4 & RDPK3SpFSAL49 & $11\,119$ & $112\,029$ & $721$ \\
\bottomrule
\end{tabular}
\caption{Timings for the Daru--Tenaud viscous shock tube with a uniform mesh.}
\label{tab:daru-timing}
\end{table}

\begin{remark}
CRKG3 remains physical with its full hyperbolic safety factor $C_i=1.0$, while CRKG4 requires the reduced factor $C_i=0.8$.
Running CRKG4 at $C_i=1.0$ leads to instability on this problem, even though its imaginary half-height $\delta(M)$ is larger than that of CRKG3 and would otherwise allow a larger hyperbolic step.
The larger steps appear too large for other components of the simulation, such as the limiter. On other problems tested, the full hyperbolic CFL remains physical even with limiting.
\end{remark}


Finally, we repeat the $\mathrm{Re}=200$ CRKG3 experiment on a graded mesh with $256\times 256$ elements at polynomial degree $N=3$.
Writing $s=(1+\xi)/2$ for the reference coordinate $\xi\in[-1,1]$, the physical map keeps $x$ uniform and uses the geometric stretch
\begin{equation}
\label{eq:daru-grade}
  \sigma(t;g)=\frac{e^{gt}-1}{e^{g}-1}
\end{equation}
\[
  y=
  \begin{cases}
    \tfrac12\,\sigma(\eta+1;g_y), & \eta\le 0,\\
    1-\tfrac12\,\sigma(1-\eta;g_y), & \eta>0,
  \end{cases}
\]
so that elements cluster toward the top and bottom walls, with grading parameter $g_y$.
Two graded meshes are considered: $g_y=1$ (mild) and $g_y=2$ (strong), along with the uniform mesh.
The corresponding smallest wall-normal element heights are $\Delta y=3.91\times10^{-3}$ (uniform), $\Delta y=2.28\times10^{-3}$ (mild), and $\Delta y=1.23\times10^{-3}$ (strong).
The CFL estimates~\eqref{eq:daru-radii}--\eqref{eq:cfl-choice} use the smallest element width in place of a uniform $h_e$.

\begin{figure}[H]
\centering
\includegraphics[width=\linewidth]{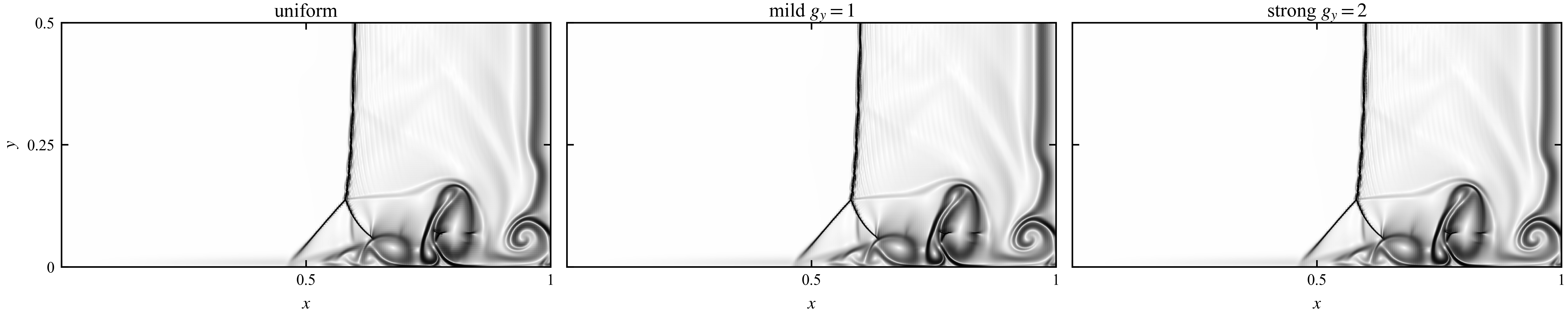}
\caption{Final numerical schlieren plots for uniformly, mildly, and strongly graded meshes at $\mathrm{Re}=200$. All solutions are qualitatively identical.}
\label{fig:daru-demos}
\end{figure}

\begin{figure}[H]
\centering
\includegraphics[width=0.92\linewidth]{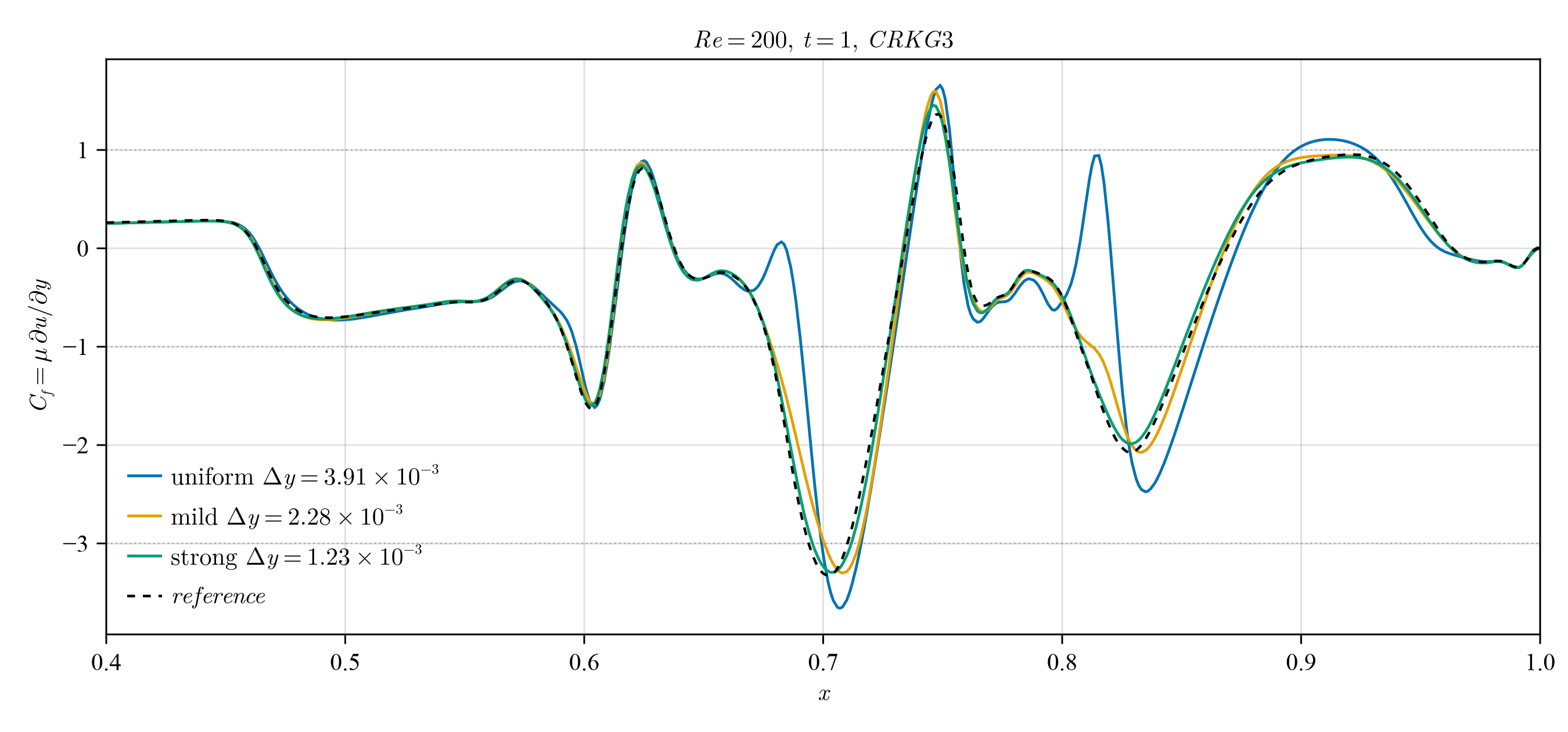}
\caption{Bottom-wall skin friction $C_f$ for uniformly, mildly, and strongly graded meshes, with the OSMP7 reference of Daru and Tenaud~\cite{tenaud2020shocktube}.}
\label{fig:daru-cf}
\end{figure}

Figure~\ref{fig:daru-demos} shows the final numerical schlieren fields on the lower half of the domain for the uniform, mild, and strong meshes; the upper half is the mirror image across $y=1/2$. The solution plots are qualitatively very similar for all three meshes.
Figure~\ref{fig:daru-cf} reports the corresponding bottom-wall skin-friction coefficient $C_f=\mu\,\partial u/\partial y$, computed with the DG derivative operator.
However, the graded meshes are significantly more accurate relative to the reference OSMP7 wall data of Daru and Tenaud~\cite{tenaud2020shocktube},
with the strongly graded mesh in particular agreeing closely with the reference data, unlike the schlieren fields, which are similar qualitatively for all meshes.
Reference data for higher $\mathrm{Re}$ is also available in~\cite{guermond2021shocktube}.

Additionally, timing is compared on the graded meshes.
The CRKG3 run is timed to completion, and then SSPRK43 is run for the same wall-clock time.
Table~\ref{tab:daru-timing-graded} reports the timings.
\begin{table}[H]
\centering
\caption{Daru--Tenaud viscous shock tube timings on graded meshes.}
\label{tab:daru-timing-graded}
\begin{tabular}{llrr}
\toprule
Mesh case & Method & Wall-clock time [s] & Simulation time \\
\midrule
mildly graded & CRKG3 & $440.23$ & $1$ \\
mildly graded & SSPRK43 & $440.26$ & $8.792\times10^{-4}$ \\
strongly graded & CRKG3 & $880.94$ & $1$ \\
strongly graded & SSPRK43 & $880.95$ & $8.319\times10^{-4}$ \\
\bottomrule
\end{tabular}
\end{table}

On the strongly graded mesh, in the same wall-clock time that CRKG3 completes the full simulation to $t=1$, SSPRK43 advances only to $t=8.319\times10^{-4}$, about $0.083\%$ of the target interval.
On the mildly graded mesh the corresponding reach is $t=8.792\times10^{-4}$.
This is still in the early stages of the boundary layer, where the parabolic stiffness is very high. The total simulation time of SSPRK43 is
expected to be faster per time unit than this early stage, but it would still be orders of magnitude slower than CRKG3.

\section{Conclusion}
\label{sec:conclusion}

We have introduced the Chain Runge--Kutta--Gegenbauer (CRKG) methods, a family of stabilized explicit Runge--Kutta schemes of orders three and four that are applicable to internal stage limiting.
The construction combines Gegenbauer chains so that the full step matches a prescribed target stability polynomial while satisfying the nonlinear order conditions, shares stages between a trunk and branches to control cost, and refines the coefficients
with a linear program and Nelder--Mead optimization for lower internal amplification.
The resulting methods have physical, monotone-in-time stages and low internal amplification, of order $10$ for CRKG3 and of order $10$--$10^{2}$ for CRKG4 over wide reaches along the negative real axis, in contrast to DUMKA3 and ROCK4, whose amplification grows to $10^{4}$ and $10^{9}$, respectively.

Numerical experiments confirm these properties.
CRKG3 and CRKG4 attain their design orders on the van der Pol oscillator.
On a one-dimensional advection--diffusion problem with stagewise positivity limiting, the methods remain physical when ROCK and DUMKA fail.
On the Daru--Tenaud viscous shock tube, CRKG is several times faster than adaptive SSPRK43 and RDPK3SpFSAL49 under the same limiter, and on a wall-graded mesh at $\mathrm{Re}=200$ the advantage over SSPRK43 is significantly more pronounced
while the wall skin friction is better resolved.


A future direction would be to apply the same chain construction to the factorized Runge--Kutta methods of O'Sullivan~\cite{osullivan2015,osullivan2017}.
We have already developed preliminary chain-based stabilized general linear methods that achieve order up to eight without splitting or extrapolation and are, to our knowledge, the highest-order stabilized methods for nonlinear problems.
These methods, however, require a significantly different construction from CRKG, including multistep constructions and complex-valued coefficients, and will be described in a separate future manuscript.

\section*{Acknowledgments}
The authors gratefully acknowledge support from the National Science Foundation under the award DMS-2607830.

\section*{AI disclosure}

During the preparation of this work, the author(s) used Claude Code and Cursor to assist with LaTeX formatting of the manuscript, assist with creating visualizations, and assistance with debugging and prototyping the numerical code. The authors have reviewed, edited, and verified these outputs and take full responsibility for the final content.



\section*{Code availability}
Implementations of the CRKG3 and CRKG4 methods, precomputed tableaus, and scripts and data for the numerical experiments in this paper are publicly available at \url{https://github.com/Sreeram-Shankar/CRKG}.

\section*{Author contributions}
\textbf{Sreeram Shankar:} Conceptualization, Methodology, Software, Validation, Formal analysis, Investigation, Writing -- Original Draft, Writing -- Review \& Editing, Visualization.

\textbf{Jesse Chan:} Conceptualization, Validation, Investigation, Resources, Writing -- Review \& Editing, Supervision, Project administration, Funding acquisition.

\appendix

\section{Mixed precision}
\label{app:mixed}
\label{sec:mixed}

The cost of a step is typically dominated by the right-hand side evaluations and by the vector operations of the recurrence, both of which are limited by memory bandwidth for large problems.
Evaluating the stages in single precision halves that traffic, and on GPUs it can also exploit faster single-precision arithmetic~\cite{abdelfattah2021,higham2022}.
A method run entirely in single precision, however, rounds the solution at every step and can lose increments smaller than the single-precision resolution, so that the solution stagnates as the step size decreases.
Thus a mixed-precision approach may be used, where the stages and right-hand sides are computed in single precision, while the solution and the accumulated update are kept in double precision~\cite{croci2023}.
Concretely, $y_n$ is rounded to single precision once at the start of the step; the stages $Y_{c,j}$ and the right-hand sides $f(Y_{c,j})$ are formed and stored in single precision; and each right-hand side is added in double precision into the update $h\sum_{c}\sum_{j}b_{c,j}f(Y_{c,j})$.
Since the update is accumulated and applied in double precision, increments are never lost to the resolution of the solution, and the solution is not rounded from step to step.

Mixed-precision time integration for explicit stabilized methods was studied by Croci and Rosilho de Souza~\cite{croci2022}, who showed that evaluating all right-hand sides in low precision limits the attainable accuracy to the order of the low-precision unit roundoff, independently of the step size, and constructed RKC methods of orders one and two that recover their full order by evaluating one or two right-hand sides in high precision.
The evaluation here is of the former kind: all right-hand sides are in single precision, so the limiting accuracy is set by the single-precision unit roundoff $u_s=2^{-24}$, scaled by the internal amplification.
Combining the order-preserving approach of~\cite{croci2022} with the CRKG construction would lead to a more robust mixed-precision evaluation. However, it becomes a more difficult problem for higher-order methods that we do not address here.

On the linear test equation, relative stage perturbations of size $u_s$ change $y_{n+1}$ by at most $\sum_j|\theta_j|\,u_s$, which is of order $\mathcal{A}\,u_s$~\cite{ketcheson2014}.
For CRKG, with $\mathcal{A}=O(1)$ to $O(10^{2})$ over the intervals of interest, this remains $O(u_s)$ to $O(10^{2}\,u_s)$, consistent with the $O(u_s)$ limit of~\cite{croci2022} and well below typical spatial discretization errors.
For DUMKA3, with $\mathcal{A}=O(10^{4})$, and for ROCK4, with $\mathcal{A}=O(10^{8})$ to $O(10^{9})$ (Section~\ref{sec:prelim-amp}), the same evaluation amplifies single-precision rounding far beyond that of CRKG, and for ROCK4 to the size of the solution itself.
Existing higher-order stabilized methods therefore cannot use this mixed-precision evaluation safely.
The same low amplification that makes stagewise limiting practical also makes mixed precision practical.
Section~\ref{sec:numerics-mixed} reports the resulting GPU speedup.

\subsection{Mixed precision comparison on a GPU}
\label{sec:numerics-mixed}

We measure the speedup of the mixed-precision evaluation of Section~\ref{sec:mixed} on an NVIDIA L40S GPU.
The test problem is the two-dimensional linear advection--diffusion equation $u_t+au_x=\nu\Delta u$ on the periodic unit square, with $a=1$ and $\nu=5\times10^{-3}$, discretized by second-order central differences on $N\times N$ grids, with exact solution $u=1+0.1\,e^{-8\pi^2\nu t}\sin(2\pi(x-at))\sin(2\pi y)$.
The right-hand side and the CRKG stage updates are implemented as fused CUDA kernels through CUDA.jl~\cite{besard2018}, with the solution and accumulated update kept in double precision and the stages and right-hand sides evaluated in single precision as in Section~\ref{sec:mixed}.
Both the double- and mixed-precision runs use the same segment count and step size, chosen once by~\eqref{eq:cfl-choice}.
After a short warm-up, one CRKG step is captured into a CUDA graph and we time $80$ graph launches on the device.
Table~\ref{tab:mp-gpu} reports the time per step, the speedup, and the errors against the exact solution.
DUMKA3 and ROCK4 are not shown here since, as discussed in Section~\ref{sec:mixed}, their amplification scales single-precision rounding to levels that are not usable.

\begin{table}[H]
\centering
\caption{Mixed-precision timing on an NVIDIA L40S.}
\label{tab:mp-gpu}
\begin{tabular}{lrrrrrrr}
\toprule
Method & $N_x$ & $M$ & $N_f$ & ms/step (f64 / mixed) & Speedup & f64 error & Mixed error \\
\midrule
CRKG3 & 2048 & 10 & 49  & 15.2 / 5.46   & 2.79 & $9.0\times10^{-8}$ & $1.4\times10^{-5}$ \\
CRKG3 & 4096 & 15 & 71  & 115.0 / 74.3  & 1.55 & $1.2\times10^{-8}$ & $2.7\times10^{-5}$ \\
CRKG3 & 8192 & 18 & 96  & 606.1 / 397.4 & 1.53 & $1.4\times10^{-9}$ & $2.9\times10^{-5}$ \\
CRKG4 & 2048 & 11 & 90  & 27.6 / 8.30   & 3.32 & $1.5\times10^{-7}$ & $5.6\times10^{-5}$ \\
CRKG4 & 4096 & 15 & 122 & 196.6 / 126.8 & 1.55 & $1.9\times10^{-8}$ & $1.3\times10^{-4}$ \\
CRKG4 & 8192 & 21 & 170 & 1072.8 / 699.8 & 1.53 & $4.2\times10^{-9}$ & $2.6\times10^{-4}$ \\
\bottomrule
\end{tabular}
\end{table}

The mixed-precision evaluation is between $1.5$ and $3.3$ times faster than the double-precision evaluation.
At $N_x=4096$ and $N_x=8192$, where the computations seem to be limited by memory bandwidth, the speedup settles at about $1.5$, consistent with the reduced memory traffic of the mixed-precision evaluation.
The larger speedup at $N_x=2048$ seems to be consistent with the working set fitting in the GPU cache.
The mixed-precision error shows the limits of the approach, of order $10^{-5}$ for CRKG3 and $10^{-4}$ for CRKG4.
It grows with the number of stages per step, in agreement with the bound $\sum_j|\theta_j|\,u_s$ of Section~\ref{sec:mixed}, and it is consistent with the limiting accuracy of order $u_s$ identified in~\cite{croci2022}.

For this smooth and well-resolved problem, the limiting accuracy exceeds the discretization error, which ranges from $10^{-7}$ to $10^{-9}$.
Mixed precision is therefore suited to simulations whose required accuracy lies above this limit, such as engineering simulations of compressible flow, and not to high-accuracy simulations of smooth solutions.

These results are preliminary and depend on the hardware and the cost of the single-precision right-hand side.
At the largest grids the speedup matches the reduction in memory traffic, which we expect to carry over to GPUs designed for scientific computing~\cite{nvidia2023l40s,nvidia2023h100}.

\section{Challenges with constructing higher-order CRKG schemes}
\label{app:crkg5}

The construction of Section~\ref{sec:crkg} extends formally beyond orders three and four.
This appendix details a potential fifth-order CRKG method and the reasons we do not pursue it.
For $p=5$, the target~\eqref{eq:target} is built from five Gegenbauer polynomials.
With the default parameters of Section~\ref{sec:crkg-target} its reach is
$\beta\approx 0.16\,L^{2}$, compared with $0.30\,L^{2}$ and $0.21\,L^{2}$ for orders three and four.
Order five adds the eight non-tall trees with five vertices to the conditions of
Table~\ref{tab:order-conditions}, so that~\eqref{eq:weight-system} carries $12$ non-tall
conditions.
Although two chains provide enough degrees of freedom in principle, designs with low
amplification and an acceptable leading error require at least six chains, and the
cost~\eqref{eq:rhs-per-step} grows accordingly.

The fifth-order target also loses stability along the imaginary axis near the origin, as described in the following proposition:
\begin{proposition}
\label{prop:imag-axis}
Let $R(z)=e^{z}+e_{p+1}z^{p+1}+e_{p+2}z^{p+2}+O(z^{p+3})$ be a polynomial of order $p$.
Then, as $y\to 0$,
\[
  |R(iy)|^{2}=
  \begin{cases}
    1+2e_{p+1}\,(-1)^{(p+1)/2}\,y^{p+1}+O(y^{p+3}), & p \text{ odd},\\[2pt]
    1+2(-1)^{p/2}\,(e_{p+1}-e_{p+2})\,y^{p+2}+O(y^{p+3}), & p \text{ even}.
  \end{cases}
\]
In particular, if $e_{p+1}<0$, then $|R(iy)|>1$ for all sufficiently small $y\neq 0$ when $p\equiv 1\pmod 4$.
\end{proposition}

\begin{proof}
With $|e^{iy}|=1$, $|R(iy)|^{2}=1+2\operatorname{Re}\bigl(e^{-iy}(e_{p+1}(iy)^{p+1}+e_{p+2}(iy)^{p+2})\bigr)+O(y^{2p+2})$.
For odd $p$, $(iy)^{p+1}=(-1)^{(p+1)/2}y^{p+1}$ is real, which gives the leading term, and the next real contributions are of order $y^{p+3}$.
For even $p$, $(iy)^{p+1}$ is imaginary, so the leading contributions come from $\operatorname{Re}(-iy\,e_{p+1}(iy)^{p+1})$ and $\operatorname{Re}(e_{p+2}(iy)^{p+2})$, both of order $y^{p+2}$.
\end{proof}

For the targets of Section~\ref{sec:crkg-target}, the leading error coefficient $e_{p+1}$ is
negative for every order, about $-8\times10^{-3}$, $-9\times10^{-4}$, and $-9\times10^{-5}$ for
$p=3$, $4$, and~$5$.
By Proposition~\ref{prop:imag-axis}, this sign makes the targets of order three and, through
$e_{5}-e_{6}<0$, order four stable near the imaginary axis, but it makes the target of order
five unstable.
The $p$ coefficients of the target are fixed by its $p$ order conditions, so the sign cannot be
corrected without additional degrees of freedom.
The instability is weak, with $|R(iy)|-1$ of order $10^{-6}$ for $|y|\le 1/2$, but it reduces the
imaginary half-height $\delta(M)$ of~\eqref{eq:delta-M} from about $3$ for order four to about
$0.15$, so that~\eqref{eq:cfl-limits} reduces the hyperbolic step size by a factor of about~$20$.

Table~\ref{tab:crkg5} compares the fifth-order target with CRKG3 and CRKG4 at the same
segment count $M=14$.
Even before counting the additional chains, the fifth-order target already has a shorter reach
per stage and an imaginary half-height about twenty times smaller than that of CRKG4.

\begin{table}[H]
\centering
\caption{Fifth-order CRKG target compared with CRKG3 and CRKG4.}
\label{tab:crkg5}
\begin{tabular}{lrrrrr}
\toprule
Method & $M$ & $\beta$ & $\beta/L^{2}$ & $\delta$ & $N_f$ \\
\midrule
CRKG3 & 14 & 520.7 & 0.30 & 2.66 & 68 \\
CRKG4 & 14 & 651.0 & 0.21 & 3.27 & 106 \\
CRKG5 & 14 & 776.1 & 0.16 & 0.15 & --- \\
\bottomrule
\end{tabular}
\end{table}

Thus a fifth-order CRKG method would have physical chain stages, but it is not practical due to the following reasons:
more chains and right-hand side evaluations per step for less reach per evaluation, larger
internal amplification, and a target unstable near the imaginary axis, which makes it
unsuitable for the problems we target.
Making that target stable near the imaginary axis would require additional degrees of freedom,
for example a sixth Gegenbauer component constrained to make $e_{6}\ge 0$, at a further cost in
reach.
For higher orders, the chain approach to the factorized Runge--Kutta methods is more promising,
as discussed in Section~\ref{sec:conclusion}.

\end{document}